\documentclass[a4paper,12pt,leqno]{article}
\usepackage{amsmath,amsfonts,latexsym}
\usepackage{amssymb}
\usepackage{graphics}

\newcommand{\nc}{\newcommand}
\nc{\tr}{{\triangle}} \nc{\vth}{{\vartheta}}
 \nc{\bt}{{\beta}}
\nc{\dl}{{\delta}} \nc{\Dl}{{\Delta}}
 \nc{\p}{{\psi}}
\nc{\gm}{{\gamma}} \nc{\Gm}{{\Gamma}} \nc{\sg}{{\sigma}}
\nc{\ve}{{\varepsilon}} \nc{\ch}{{\cal H}} \nc{\cf}{{\cal F}}
\nc{\cp}{{\cal P}}
 \nc{\td}{\tilde}

\newtheorem{rem}{Remark}[section]

\newtheorem{cor}{Corollary}[section]
\newtheorem{Lemma}{Lemma}[section]

\begin{document}

%***************************************************************

%***************************************************************************

\title
{Recursive Parameter Estimation: \\Asymptotic expansion.}
   \author{Teo Sharia}
\date{}
\maketitle
\begin{center}
{\it
Department of Mathematics \\Royal Holloway,  University of London\\
Egham, Surrey TW20 0EX \\ e-mail: t.sharia@rhul.ac.uk }
\end{center}
\begin{abstract}
 We consider
 estimation procedures which are recursive in the sense that each successive
 estimator is obtained from the previous one by a simple adjustment.
The model considered in the paper is very general as we do not
impose any  preliminary  restrictions on the probabilistic nature of
the observation process and cover a wide class of nonlinear
recursive procedures. In this paper we study asymptotic behaviour of the
recursive estimators.
The results of the paper can be used to determine
the form of  a recursive procedure  which is expected to have the same asymptotic
properties as the corresponding non-recursive one defined as a
solution of the corresponding estimating equation.
\end{abstract}
\begin{center}
\end{center}

\section{Introduction}

Let $X_1, \dots, X_n$ be independent identically distributed
(i.i.d.) random variables (r.v.'s) with  a common distribution
function
  $F_{\theta}$ with a  real unknown parameter $\theta$. An
$M$-estimator of $\theta$ is defined as a statistic $
\hat\theta_n=\hat\theta_n(X_1,\dots,X_n), $ which is a solution
w.r.t. $v$ of the estimating equation
%%%%%%%%%%%%%%%%%%%%%%%%%%%%%%%%%%%%%%%%%%%%%%%%%
%                       (esteq)
%%%%%%%%%%%%%%%%%%%%%%%%%%%%%%%%%%%%%%%%%%%%%%
\begin{equation} \label{esteq}
\sum_{i=1}^n \psi(X_i;v)=0,
\end{equation}
where $\psi$ is a suitably chosen function. For example, if
$\theta$ is a location parameter in  the normal family of
distribution functions, the choice   $\psi(x,v)=x-v$ gives
 the MLE (maximum likelihood estimator).
For the same problem, if $\psi(x,v)=\mbox{sign}(x-v),$ the
solution of \eqref{esteq} reduces to the median of $X_1, \dots,
X_n$. In general, if $f(x,\theta)$ is the probability density
function (or probability function) of  $F_{\theta}(x)$  (w.r.t. a
$\sigma$-finite measure $\mu$)  then the choice
$\psi(x,v)=f'(x,v)/f(x,v)$ yields the MLE.

 Suppose now that $X_1, \dots, X_n$ are not
necessarily independent or identically  distributed r.v's, with a
joint distribution depending on a real parameter $\theta$. Then an
$M$-estimator of $\theta$ is defined as  a solution of the
estimating equation
%%%%%%%%%%%%%%%%%%%%%%%%%%%%%%%%%%%%%%%%%%%%%%%%%%%%
%                  \eqref{esteqg}
%%%%%%%%%%%%%%%%%%%%%%%%%%%%%%%%%%%%%%%%%%%%%%%%%%%%%%
\begin{equation}\label{esteqg}
\sum_{i=1}^n \psi_i(v)=0,
\end{equation}
where   $\psi_i(v)=\psi_i(X_{i-k}^i; v)$ with $X_{i-k}^i=
(X_{i-k}, \dots,X_i )$. So, the $\psi$-functions  may now depend
on the past observations as well. For instance, if $X_i$'s are
observations from a discrete time Markov process, then one can
assume that   $k=1$. In general, if no restrictions are placed on
the dependence structure of the process $X_i$, one may need to
consider $\psi$-functions depending on the vector of all past and
present observations of the process (that is, $k=i-1$). If the
conditional probability density function (or probability function)
of the observation $X_i,$ given $X_{i-k}, \dots,X_{i-1},$ is
$f_i(x,\theta)=f_i(x,\theta|X_{i-k}, \dots,X_{i-1})$,  then one
can obtain the  MLE  on choosing
$\psi_i(v)=f'_i(X_i,v)/f_i(X_i,v).$ Besides MLEs, the class of
$M$-estimators includes estimators
 with special properties such as {robustness}.
 %An estimator
 %is  robust if its behaviour is not ``seriously
 %affected'' by violations of underlying assumptions.
Under certain regularity and ergodicity conditions, it can be
proved that there exists a consistent sequence of solutions of
\eqref{esteqg} which has  the property of local asymptotic
linearity. (A comprehensive bibliography can be found in, e.g.,
Hampel at al (1986) and Rieder (1994).)

 If   $\psi$-functions are
 nonlinear, it is rather difficult to
 work with  the corresponding estimating equations, especially
  if for every sample size $n$
 (when new data are acquired),
 an estimator has to be computed afresh. In this paper we consider
 estimation procedures which are recursive in the sense that each successive
 estimator is obtained from the previous one by a simple adjustment.
 Note that for a linear estimator,
 e.g., for the sample mean, $\hat\theta_n=\bar X_n$ we have
$\bar X_n=(n-1)\bar X_{n-1}/n+X_n/n$, that is
$\hat\theta_n=\hat\theta_{n-1}(n-1)/n+X_n/n$, indicating that the
estimator $\hat\theta_n$ at  each  step $n$ can be obtained
recursively using the estimator at the previous step
$\hat\theta_{n-1}$ and the new information $X_n$. Such an exact
recursive relation may not hold for nonlinear estimators (see,
e.g., the case of the median).

In general, the following heuristic argument can be used to
establish   a possible form of an approximate recursive relation
(see also Jure$\check{{\rm c}}$kov$\acute{{\rm a}}$ and Sen (1996),
 Khas'minskii  and Nevelson (1972),
 Lazrieva and  Toronjadze (1987)).
 Since  $\hat\theta_n$ is defined
as   a root  of the estimating equation \eqref{esteqg}, denoting
the left hand side of \eqref{esteqg} by $M_n(v)$ we have
$M_n(\hat\theta_n)=0$ and  $M_{n-1}(\hat\theta_{n-1})=0$. Assuming
that the difference $\hat\theta_n-\hat\theta_{n-1}$ is ``small''
we can write
$$
0=M_n(\hat\theta_n)-M_{n-1}(\hat\theta_{n-1})=
M_n\left(\hat\theta_{n-1}+(\hat\theta_n-\hat\theta_{n-1})\right)-
M_{n-1}(\hat\theta_{n-1})
$$
$$
\approx
M_n(\hat\theta_{n-1})+M_n'(\hat\theta_{n-1})(\hat\theta_n-\hat\theta_{n-1})
-M_{n-1}(\hat\theta_{n-1})
$$
$$
=M_n'(\hat\theta_{n-1})(\hat\theta_n-\hat\theta_{n-1})
+\psi_n(\hat\theta_{n-1}).
$$
Therefore,
$$
\hat\theta_n\approx\hat\theta_{n-1}-\frac{\psi_n(\hat\theta_{n-1})}
{M_n'(\hat\theta_{n-1})},
$$
where $M_n'(\theta)=\sum_{i=1}^n \psi'_i(\theta)$. Now, depending
 on the nature of the underlying model,  $M_n'(\theta)$ can be replaced by
 a simpler expression. For instance, in  i.i.d.
 models with $\psi(x,v)=f'(x,v)/f(x,v)$  (the  MLE case), by the strong law of large
 numbers,
 $$
 \frac {M_n'(\theta)}n =\frac1n  \sum_{i=1}^n
 \left(f'(X_i,\theta)/f(X_i,\theta)\right)'\approx E_\theta
 \left[  \left(f'(X_1,\theta)/f(X_1,\theta)\right)'\right]
 =-i(\theta)
 $$
 for large $n$'s,
where $i(\theta)$ is the one-step Fisher   information. So, in
this case, one can use the recursion\footnote{This procedure should not be confused with the
Newton-Raphson iterative method. See the corresponding discussion in the Introduction of Sharia (2006a).}

%%%%%%%%%%%%%%%%%%%%%%%%%%%%%%%%%%%%%%%%%%%%%%%%%%
%                 {mliid}
%%%%%%%%%%%%%%%%%%%%%%%%%%%%%%%%%%%%%%%%%%%%%%%%%
\begin{equation} \label{mliid}
\hat\theta_n=\hat\theta_{n-1}+\frac{1}{n \; i(\hat\theta_{n-1})}
\frac{f'(X_n,\hat\theta_{n-1})}{f(X_n,\hat\theta_{n-1})},
~~~~~~~~~n\ge 1,
\end{equation}
to construct an estimator which is ``asymptotically equivalent''
to the MLE.

Motivated by  the above argument,   we consider  a class of
estimators
%%%%%%%%%%%%%%%%%%%%%%%%%%%%%%%%%%%%%%%%%%%%%%%%%%
%                   \eqref{rec1}
%%%%%%%%%%%%%%%%%%%%%%%%%%%%%%%%%%%%%%%%%%%%%%%%%
\begin{equation} \label{rec1}
\hat\theta_n=\hat\theta_{n-1}+
{\Gamma_n^{-1}(\hat\theta_{n-1})}\psi_n(\hat\theta_{n-1}),
~~~~~~~~~n\ge 1,
\end{equation}
where $\psi_n$ is a suitably chosen vector process, $\Gamma_n$  is
a (possibly random) normalizing matrix process and
$\hat\theta_0\in {\mathbb{R}}^m$ is some initial value.
 If the
conditional probability density function (or the probability
function) of the observation $X_n,$ given $X_{1}, \dots,X_{n-1},$
is $ f_n(\theta, x|x_1^{n-1})=f_n(x,\theta|x_{1}, \dots,x_{n-1})$,
then one can obtain a   ML (maximum likelihood) type recursive
estimator  on choosing $\psi_n(\theta)=\dot f^T_n(\theta,
X_n|X_1^{n-1})/f_n(\theta, X_n|X_1^{n-1})$ (the dot denotes the
row-vector of partial derivatives w.r.t. $\theta \in {\mathbb{R}}^m$ and $T$ is the transposition).

 Note that
while the main goal is  to study recursive procedures with
non-linear $\psi_n$ functions, it is worth mentioning that any
linear estimator can be written in  the form  \eqref{rec1} with
linear, w.r.t. $\theta$, $\psi_n$ functions. Indeed, if
$\hat\theta_n=\Gamma_n^{-1}\sum_{k=1}^n h_{k} (X_k),$  where
$\Gamma_k$ and $h_k(X_k)$ are  matrix and vector processes  of
suitable dimensions, then (see Section 4.2 for details)
$$
\hat\theta_n=\hat\theta_{n-1}+\Gamma_n^{-1}\left( h_n(X_n)-
(\Gamma_n-\Gamma_{n-1})\hat\theta_{n-1}\right),
$$
which is obviously of the form  \eqref{rec1} with $\psi_n(\theta)=
h_n(X_n)- (\Gamma_n-\Gamma_{n-1})\theta.$

 Note also
that in the iid case,
 \eqref{mliid} can be regarded as a stochastic iterative scheme,
i.e.,  a classical stochastic approximation procedure,  to detect
the root of an unknown function when the latter can only be
observed with random errors (see Remark 3.1 in Sharia (2006a)). A theoretical
implication of this is that by studying the procedures
\eqref{mliid}, or in general \eqref{rec1}, we study  asymptotic
behaviour of the estimator of the unknown parameter. As far as
applications  are concerned, there are several advantages in using
\eqref{rec1}. Firstly, these procedures are easy to use since each
successive
 estimator is obtained from the previous one by a simple adjustment
and without  storing all the data unnecessarily. This is
especially convenient when the data come sequentially.
 Another potential benefit of using \eqref{rec1}  is that it allows one to
 monitor and detect certain changes  in probabilistic characteristics  of
 the underlying process such as change of the value of the unknown parameter. So, there
 may be a benefit in using these procedures in linear cases as well.

In  i.i.d. models,
 estimating procedures similar to \eqref{rec1} have been studied by a number of
authors using methods of stochastic approximation theory (see,
e.g., Khas'minskii  and Nevelson (1972), Fabian (1978),
 Ljung and Soderstrom (1987), Ljung et al (1992),  and references therein). Some work
has been done for  non i.i.d. models as well. In particular,
Englund et al (1989) give an asymptotic representation results for
certain type of  $X_n$ processes. In Sharia (1998), theoretical
results on convergence, rate of convergence and the asymptotic
representation are given under certain regularity and ergodicity
assumptions on the model, in the one-dimensional  case with
$\psi_n(x,\theta)=\frac{\partial}{\partial\theta} \mbox{log}
f_n(x,\theta)$ (see also Campbell (1982), Sharia (1992),
  and  Lazrieva et al (1997)).

We study multidimensional  estimation procedures of type
\eqref{rec1} for the general statistical model. In  Sharia (2006a),
 imposing ``global'' restrictions on the processes $\psi$
and $\Gamma$, we study ``global'' convergence of the recursive
estimators, that is the convergence for an arbitrary starting
value $\hat\theta_0$. In Sharia (2006b),   we present results on
the rate of the convergence. In this paper we are concerned with
asymptotic behaviour  of the estimators defined by \eqref{rec1}.
Since the model considered is very general, the main objective is
to prove that $\hat \theta_n$ is locally asymptotically linear,
that is, for each $\theta$ there exist a matrix process
$G_n(\theta)$ such that
$$
\hat \theta_n -\theta =G^{-1}_n(\theta) \sum_{i=1}^n
\psi_i(\theta) + \varepsilon_n^\theta,
$$
where $ G^{1/2}_n(\theta) \varepsilon_n^\theta \to 0 $ in
probability $P^\theta$ (see Section 2 for a more general
definition).

Since $\psi_t(\theta)$ is typically a martingale-difference,
  asymptotic distribution  of an asymptotically  linear estimator can be studied
using a suitable form  of the central limit theorem for
martingales (see e.g., Feigin (1985), Hutton and Nelson (1986),
Jacod and Shiryayev (1987). Detailed discussion of the literature
on this subject can be found in
 Barndorff-Nielsen and  Sorensen (1994),  Heyde (1997) and Prakasa-Rao  (1999)).
For example, results in Shiryayev (1984) (see, e.g.,  Ch.VII,
$\S$8, Theorem 4)  show that under certain conditions, local
asymptotic linearity  implies asymptotic normality. In the
standard case of i.i.d. observations, assuming that
$$
\p_n(\theta)=\p(\theta,X_n)
$$
 has zero mean and a finite
second moment and   $G_n(\theta)=n\gamma(\theta),$ for some
non-random  invertible $\gamma(\theta)$, it follows that
 $$
{\cal L}\left( n^{1/2}(\hat  \theta_n -\theta) \mid P^\theta
 \right) {\stackrel{w}{\to}} {\cal N} \left(0,\;\; \gm^{-1}(\theta)
 j_\p(\theta) \gm^{-1}(\theta)\right)
$$
where
$$
   j_\p(\theta)=\int
\p(\theta,x)\p^T(\theta,x) f(\theta,x)\mu(\,dx)< \infty.
$$
In particular, in the  case of likelihood recursion with
$$
\psi(\theta,x)=\dot f^T(\theta,x)/f(\theta,x),
$$
if $\gamma(\theta)$ is the one-step Fisher information, that is,
$$
\gamma(\theta)=i(\theta)=j_{\dot f^T/f}(\theta),
$$
it follows  that  $\hat \theta_n$ is asymptotically normal with
parameters  $ (0, i^{-1}(\theta))$, i.e.
$$
{\cal L}\left( n^{1/2}(\hat  \theta_n -\theta) \mid P^\theta
 \right) {\stackrel{w}{\to}} {\cal N} \left(0,\;\; i^{-1}(\theta)\right),
$$
 meaning that  $\hat \theta_n$  is
asymptotically efficient.
%in the sense of having asymptotically
%minimum variance among the class of all asymptotic normal
%estimators.
In general, in the case of one dimensional parameter $\theta$, an
estimator is said to be {\it asymptotically efficient} if it is
asymptotically linear with
$$
\p_n(\theta)= \dot f_n(\theta, X_n|X_1^{n-1})/f_n(\theta,
X_n|X_1^{n-1}) \;\;\;  \mbox{and} \;\;\; G_n(\theta)=I_n(\theta).
$$
where $ I_n(\theta) $ is the conditional Fisher information. This
kind of efficiency is called asymptotic first order efficiency.
The motivation behind this general definition is the same as in
the classical scheme  of i.i.d. observations. For a detailed
discussion of this notion see, e.g.,
  Hall and  Heyde (1980), Section 6.2. Under relatively mild conditions,
asymptotically  efficient estimators are asymptotically equivalent to
the MLE $T_n$, i.e.
$$
I_n^{1/2}(\theta)(\hat\theta_n-T_n) \to 0
$$
in probability (see, e.g., Hall and Heyde (1980), Section 6.2,
Theorem 6.2.). For the generalisation of these concepts see Heyde
(1997).

It is worth mentioning that the global convergence results for
\eqref{rec1} were obtained in Sharia (2006a) under conditions that
allow $\Gamma_n$ to belong to quite a wide class of processes
which does not directly depend  on the choice of $\p_n$'s (see Remark
3.1 below). In
order to study the rate of convergence, one has to restrict the
class of allowed $\Gamma_n$'s (see Sharia (2006b)). It turns out
that when dealing with local asymptotic linearity, one has to
restrict this class even further - to an explicit  choice of
$\Gamma_n$, depending on the choice of $\p_n$ (see Remark
3.2(iv)--(vii) below). In other words, the  results of the paper
tell one how to construct a locally asymptotically linear
procedure \eqref{rec1} with given $\p_n$'s. The fact that one is
restricted to this choice of $\Gamma_t$ is probably not very
surprising in retrospective, but this issue does not seem to have
been discussed in the existing literature.

An estimator defined by \eqref{rec1} is a  recursive analogue of
the corresponding $M$-estimator   defined as a solution of the
estimating equation \eqref{esteqg}. It should also be noted that
the recursive procedure \eqref{rec1} is not a numerical solution
of \eqref{esteqg}. Nevertheless, under quite mild conditions, the
recursive estimator and the corresponding $M$-estimator are
expected to have the same (or equivalent) asymptotic linearity
expansions. It therefore follows that they are  asymptotically
equivalent, in the sense that, depending on the regularity and
ergodicity properties of the underlying model, they both have the
same asymptotic distribution.

The paper is organized as follows.  Section 2 introduces the main
objects and definitions. The main results are obtained in Section
3 with various comments and explanations of the conditions used
there.  In Section 4 we give examples to illustrate the results of
the paper.

%\newpage

                       %S E C T I O N 2

\section {Basic  model}
\setcounter{equation}{0}

Let   $X_t,\;\; t=1,2,\dots ,$
 be  observations taking values in a  measurable space
$({\bf X},{\cal B}({\bf X}))$  equipped with  a $\sigma$-finite
measure $\mu.$ Suppose  that the distribution of the process $X_t$
depends on an unknown parameter $\theta \in \Theta,$ where
$\Theta$ is an open subset of the $m$-dimensional Euclidean space
$\mathbb{R}^m$. Suppose also  that for each  $t=1,2,\dots$, there
exists a  regular conditional probability density of $X_t$ given
values of past observations  of $X_{t-1},\dots , X_2, X_1$, which
will be denoted by
$$
f_ t(\theta, x_t \mid x_1^{t-1})= f_ t(\theta, x_t \mid
x_{t-1},\dots ,x_1),
$$
where $ f_ 1(\theta, x_1 \mid x_1^0)= f_ 1(\theta, x_1) $ is the
probability  density  of the random variable $X_1.$ Without loss
of generality we assume that all random variables are defined on a
probability space $
 (\Omega, \cf)
 $
and denote by $ \left\{P^\theta, \; \theta\in \Theta\right\} $ the
family of the corresponding distributions on $
 (\Omega, \cf).
 $

  Let  $\cf_t=\sigma (X_1,\dots ,X_t)$ be the $\sigma$-field
generated by the random variables $ X_1,\dots ,X_t.$ By
$\left(\mathbb{R}^m, {\cal B}( \mathbb{R}^m ) \right)$ we denote
the $m$-dimensional Euclidean space with the Borel
$\sigma$-algebra ${\cal B}( \mathbb{R}^m )$. Transposition of
matrices and vectors is denoted by $T$. By $(u,v)$ we denote the
standard scalar product of $u,v \in \mathbb{R}^m,$ that is,
$(u,v)=u^Tv,$ and the corresponding norm is denoted by $\|u\|$.

Suppose that $h$ is a real valued function  defined on
 $\Theta\subset  {{\mathbb{R}}}^m$.  We denote by $\dot h(\theta)$  the row-vector
 of partial derivatives
of $h(\theta)$ with respect to the components of $\theta$, that
is,
 $$
 \dot h(\theta)=\left(\frac{{\partial}}{{\partial} \theta^1} h(\theta), \dots,
 \frac{{\partial}}{{\partial} \theta^m} h(\theta)\right).
 $$
%Also we denote by  \"{h}$(\theta)$ the  matrix of second partial derivatives.
 The $m\times m$ identity matrix is denoted by ${{\bf 1}}$.

 If for each  $t=1,2,\dots$, the derivative
$ \dot f_t(\theta, x_t \mid x_1^{t-1})$ w.r.t. $\theta$
  exists, then we
can define
$$
l_t(\theta, x_t \mid x_1^{t-1})=\frac 1 {f_ t(\theta, x_t \mid
x_1^{t-1})} \dot f_t^T(\theta, x_t \mid x_1^{t-1})
$$
 and the process
$$
l_t(\theta)=l_t(\theta, X_t \mid X_1^{t-1})
$$
(with the convention $0/0=0$). Let us denote
$$
i_t(\theta\mid x_1^{t-1}) =\int l_t(\theta,z\mid x_1^{t-1})
l^{T}_t(\theta,z\mid x_1^{t-1})f_t(\theta,z\mid x_1^{t-1})\mu
(dz).
$$
The {\it one step   conditional  Fisher information matrix} for
$t=1,2,\dots$ is defined as
$$
i_t(\theta)=i_t(\theta\mid X_1^{t-1}).
$$
Note that  the process $i_t(\theta)$  is {\it ``predictable''},
that is, the random variable $i_t(\theta),$ is $\cf_{t-1}$
measurable for each $t\ge 1.$
 Note also that by   definition,
$i_t(\theta)$ is   a version of the conditional expectation w.r.t.
${\cal{F}}_{t-1},$  that is,
$$
i_t(\theta)= E_\theta\left\{l_t(\theta) l^{T}_t(\theta) \mid
{\cal{F}}_{t-1}\right\}.
$$
Everywhere in the present work conditional expectations  are meant
to be
 calculated as integrals w.r.t. the conditional probability densities.

The {\it conditional  Fisher information} at time $t$ is
$$
I_t(\theta)=\sum_{s=1}^t i_s(\theta),  \;\;\;\;\;\;\;\;\;
t=1,2,\dots.
$$

We say that  ${\bf \psi}= \{\psi_t(\theta, x_t, x_{t-1}, \dots,
x_1)\}_{t\ge 1}$ is a sequence of estimating functions and write
$\bf \psi \in \Psi$, if for each ${t\ge 1},$ $ \psi_t(\theta, x_t,
x_{t-1}, \dots, x_1)  :
 \Theta \times {\bf X}^t \;\;\to \;\;{\mathbb{R}}^m
$ is a   Borel function.

Let $\psi \in {\bf\Psi}$ and denote $\psi_t(\theta)=\psi_t(\theta,
X_t, X_{t-1},\dots, X_1).$ We write $\bf \psi \in \Psi^M$ if  ~
$\psi_t(\theta)$ is a martingale-difference process for each
$\theta \in \Theta,$ ~ i.e.,~ if $
E_\theta\left\{\psi_t(\theta)\mid {\cal{F}}_{t-1}\right\}=0$ for
each $t=1,2,\dots $ (we assume that the conditional expectations
above are well-defined and ${\cal{F}}_{0}$ is the trivial
$\sigma$-algebra).

Note   that if differentiation of the equation
$
1=\int  f_t(\theta, z \mid x_1^{t-1}) \mu(dz)
$
is allowed under the integral sign, then $ \{l_t(\theta)\}_{t\ge 1}\in {\bf\Psi^M}$.

Suppose that  $\psi \in {\bf\Psi}$  and $\Gamma_t(\theta)$ is
 a predictable $m\times m$ matrix process (i.e. a matrix with predictable components  $\Gamma_t^{ij}(\theta)$ )
 with  $\mbox{det} \Gamma_t(\theta)\neq 0.$
We say that an estimator $\hat\theta_t$ is
 {\it locally asymptotically  linear} if
for each $\theta \in \Theta, $
\begin{equation}\label{LinSt}
\hat \theta_t =\theta + \Gm^{-1}_t(\theta) \sum_{s=1}^t
\p_s(\theta)+\ve_t^\theta,
\end{equation}
and $A_t(\theta) \ve_t^\theta \to 0$  in probability  $P_\theta,$ where
$A_t(\theta)$ is a sequence of  $m\times m$ matrices    such that
$ A_t(\theta)\to \infty $  in probability $P^\theta,$
 and $A_t(\theta)
\Gm_t^{-1}(\theta)A_t(\theta)\to {\eta(\theta)} $ weakly w.r.t.
$P^\theta$ for some random matrix $\eta(\theta).$
That is, $\hat\theta_t$ is
  locally asymptotically  linear if
\begin{equation}\label{aslin}
A_t(\theta)(\hat \theta_t^* -\hat\theta_t)\to 0
\end{equation}
in probability $P^\theta$, where
\begin{equation}\label{LinSt}
\hat \theta_t^* =\theta + \Gm^{-1}_t(\theta) \sum_{s=1}^t
\p_s(\theta) ,
\end{equation}
is a linear statistic.

%If $\bf \psi \in \Psi^M$,  asymptotic behaviour of  linear, and
%therefore, locally asymptotically linear estimators can be studied
%using a suitable form of the central limit theorem for martingales.

\vskip+0.5cm

 {\bf Convention} {\it  Everywhere in the present work
 $\theta\in \mathbb{R}^m $ is an arbitrary but fixed value
of the parameter.  Convergence and all relations between random
variables are meant with probability one w.r.t. the measure
$P^\theta$ unless specified otherwise. A sequence of random
variables $(\xi_t)_{t\ge1}$ has some property eventually if for
every $\omega$ in a set $\Omega^\theta$ of $P^\theta$ probability
1,
 $\xi_t$  has this property for all $t$ greater than some
$t_0(\omega)<\infty$.}

%\newpage

  %           3   Main results

\section{Main results}
\setcounter{equation}{0} Suppose   that $\bf \psi \in \Psi$ and
$\Gm_t(\theta)$,  for each $\theta\in \mathbb{R}^m$, is  a
predictable $m\times m$ matrix process with $ \mbox{det}
~\Gm_t(\theta)\neq0$, $t\ge 1$.
 Consider  the estimator $\hat \theta_t$  defined by
                       %(rec2)
\begin{equation}\label{rec2}
\hat \theta_t=\hat \theta_{t-1}+\Gm_t^{-1}(\hat
\theta_{t-1})\p_t(\hat \theta_{t-1}), \qquad
 t\ge 1,
\end{equation}
where  $\hat \theta_0\in\mathbb{R}^m $ is an arbitrary  initial
point.

Let $\theta\in \mathbb{R}^m $ be an arbitrary but fixed value of
the parameter  and  for any  $u\in \mathbb{R}^m$ define
$$
R_t(\theta,u)=\Gm_t(\theta)\Gm_t^{-1}(\theta+u)E_\theta
\left\{\psi_t(\theta+u)\mid {\cf}_{t-1}\right\}.
$$
 Denote  $\Dl_t=\hat\theta_t-\theta$. Then \eqref{rec2} can
be rewritten as
                     %Dl
\begin{equation}\label{Dl}
\Dl_t=\Dl_{t-1}+\Gm^{-1}_t(\theta) R_t(\theta,\Dl_{t-1}) +
\Gm^{-1}_t(\theta)\ve_{\theta t},
\end{equation}
where
$$
\ve_{\theta t}=\Gm_t(\theta)\Gm_t^{-1}(\theta+\Dl_{t-1})\p_t(\theta+\Dl_{t-1})-
R_t(\theta,\Dl_{t-1})
$$
is a $P^{\theta}$-martingale difference.

  Let $\Dl_0^*=0$ and for $t\ge 1$ denote
  $\Dl_t^*=\hat\theta_t^*-\theta$ where  $\hat\theta_t^*$ is defined by \eqref{LinSt}.
Then,
\begin{eqnarray}\label{solut}
\Dl_t^*-\Dl_{t-1}^*
& &=\Gm^{-1}_t(\theta) \sum_{s=1}^t \p_s(\theta)
- \Gm^{-1}_{t-1}(\theta) \sum_{s=1}^{t-1}\p_s(\theta) \nonumber \\
& & =\left (\Gm^{-1}_t(\theta)-\Gm^{-1}_{t-1}(\theta)\right) \sum_{s=1}^{t-1}\p_s(\theta)
+ \Gm^{-1}_t(\theta)\p_t(\theta) \\
& & =\Gm^{-1}_t(\theta)\left (\Gm_{t-1}(\theta)-\Gm_{t}(\theta)\right)\Dl_{t-1}^*
+ \Gm^{-1}_t(\theta)\p_t(\theta). \nonumber
\end{eqnarray}
It therefore follows   that
$\Dl_t^*$ satisfies the recursive relation given by
                                  %(Dl*)
\begin{equation}\label{Dl*}
\Dl _t^*=\Dl_{t-1}^*- \Gm_t^{-1}(\theta)\tr\Gm_t(\theta) \Dl
_{t-1}^*+\Gm_t^{-1}(\theta)\ve_{\theta t}^*, ~~~~~~~~~ t\ge 1,
\end{equation}
where  $\tr\Gm_t(\theta)=\Gm_{t}(\theta)-\Gm_{t-1}(\theta)$ and $\ve_{\theta t}^*=\p_t(\theta)$. By comparing
equations \eqref{Dl} and \eqref{Dl*}, one can obtain the following
result on the asymptotic relationship between $\hat\theta_t$ and
$\hat\theta_t^*.$

%%%%%%%%%%%%%%%%%%%%%%%%%%%%%%%%%%%%%%%%%%%%%%%%%%%%%%%%%%%%%%%%%%%%
                    % L e m m a  3.1
%%%%%%%%%%%%%%%%%%%%%%%%%%%%%%%%%%%%%%%%%%%%%%%%%%%%%%%%%%%%%%%%%%%%
\begin{Lemma}
Suppose that $\bf \psi \in \Psi$ and there exists
a sequence of   invertible random  matrices $A_t(\theta)$ such
that $A_t^{-1}(\theta)\to 0$  in probability $P^\theta$ and
\begin{description}
\item[(E)]
 $$
A_t(\theta) \Gm_t^{-1}(\theta)A_t(\theta)\to \eta(\theta)
$$
weakly w.r.t. $P^\theta,$ where $\eta(\theta)$ is a random matrix
with $\eta(\theta) < \infty $ $P^\theta$-a.s.;
 \item[(1)]
 $$
\lim_{t\to\infty}
 A_t^{-1}(\theta)  \sum_{s=1}^t \left(\tr\Gm_s(\theta)\Dl_{s-1}
 +R_s(\theta,\Dl_{s-1})
\right) = 0
$$
in probability $P^\theta$;
\item[(2)]
 $$
  \lim_{t\to\infty} A_t^{-1}(\theta)\sum_{s=1}^t {\cal E}_s(\theta) = 0
$$
in probability $P^\theta$, where
$$
{\cal E}_s(\theta)=\Gm_s(\theta)\Gm^{-1}_s(\theta+\Dl_{s-1})
\left\{\p_s(\theta+\Dl_{s-1})-E_\theta
\left\{\psi_s(\theta+\Dl_{s-1})\mid {\cf}_{s-1}\right\})\right\}
-\p_s(\theta).
$$
\end{description}
Then $ A_t(\theta)(\hat \theta_t^* -\hat\theta_t)\to 0 $ in
probability $P^\theta$ (i.e., $\hat \theta_t^*$ is locally asymptotically linear).
\end{Lemma}
%%%%%%%%%%%%%%%%%%%%%%%%%%%%%%%%%%%%%%%%%%%%%%%%%%%%%%%%%%%%%%%%%%%%
                  % P R O O F     of   L e m m a  3.1
%%%%%%%%%%%%%%%%%%%%%%%%%%%%%%%%%%%%%%%%%%%%%%%%%%%%%%%%%%%%%%%%%%%%
{\bf Proof.} To simplify notation we drop the fixed  argument or
the index $\theta$ in some of the expressions below. Denote
$ \delta_t:= \hat \theta_t-\hat \theta_t^*=\Dl_t-\Dl_t^*. $
 Subtraction  \eqref{Dl*} from \eqref{Dl} yields the recursive relation
%                             (dif)
\begin{equation}\label{dif}
\delta_t=\left({\bf 1}- \Gm_t^{-1}\tr\Gm_t\right)\delta_{t-1}+
 \Gm_t^{-1}
(\ve_t-\ve_t^*)+
 \Gm_t^{-1}(\tr\Gm_t\Dl_{t-1}+R_t(\theta,\Dl_{t-1})) .
\end{equation}
Denote $\ch_t:=\sum_{s=1}^t\;\left(\tr\Gm_s(\theta)\Dl_{s-1}
 +R_s(\theta,\Dl_{s-1})
\right)$ and  $M_t:=\sum_{s=1}^t \; [\ve_s-\ve_s^*].$ Then the
expression
$$
\delta_t=\Gm_t^{-1}\left\{M_t+ \ch_t+\delta_0\right\},
~~~~~~~~t\ge 1
$$
can easily be obtained
 by inspecting the difference between $t$'th
and $(t-1)$'th term of this sequence (exactly in the same way as in \eqref{solut}), to check that \eqref{dif} holds.

 Now, (1)  implies that $ A_t^{-1}\ch_t\to 0 $ in
probability $P^\theta$. Also, by (2),
$A_t^{-1}M_t=A_t^{-1}(\theta)\sum_{s=1}^t {\cal E}_s(\theta) \to
0$ in probability $P^\theta$. So, using (E), it follows that   $A_t\dl_t \to 0$ in probability
$P^\theta$.
 $\diamondsuit $

 \bigskip

Next result gives  sufficient conditions for (1) and
(2).

\medskip
\noindent
 {\bf Proposition  3.1}

 \noindent
{\em {\bf (a)}
Suppose  that $A_t(\theta)$  in Lemma 3.1  are diagonal matrices  with  non-decreasing (w.r.t. $t$) elements   and
\begin{description}
\item[(L1)]
$$
A_t^{-2}(\theta)\sum_{s=1}^t
A_s(\theta)[\tr\Gm_s(\theta)\Dl_{s-1}+R_s(\theta,\Dl_{s-1})] \to 0
$$
in probability $P^\theta$;
\end{description}
Then (1) holds.

\medskip

\noindent
{\bf (b)}
Suppose  that $A_t(\theta)$  in Lemma 3.1  are diagonal   non-random matrices,
 $\bf \psi \in \Psi^M$ and
\begin{description}
\item[(L2)]
$$
 \lim_{t\to\infty} \frac1{(A_t^{(jj)}(\theta))^2}  \sum_{s=1}^t
E_\theta \left\{ \left({\cal E}_s^{(j)}(\theta)\right)^2
\mid{\cf}_{s-1}\right\} = 0
$$
in probability $P^\theta$, where $A_t^{(jj)}(\theta)$ is the $j$-th diagonal element of the matrix  $A_t(\theta)$
and ${\cal E}_s^{(j)}(\theta)$ is the $j$-th component of
$
{\cal E}_s(\theta)$ which is defined in (2).
\end{description}
Then (2) holds.

\medskip

\noindent
{\bf (c)}
Suppose  that $A_t(\theta)$  in Lemma 3.1  are diagonal with  non-decreasing elements
$A_t^{(jj)}(\theta) \to \infty,$   $\bf \psi \in \Psi^M$   and
\begin{description}
\item[(LL2)]
$$
 \sum_{s=1}^\infty
\frac {E_\theta \left\{ ({\cal E}_s^{(j)}(\theta))^2
\mid {\cf}_{s-1}\right\}} {(A_s^{(jj)}(\theta))^2} < \infty
$$
$P^\theta$-a.s., where ${\cal E}_s^{(j)}(\theta)$ is the $j$-th component of
$
{\cal E}_s(\theta)$ which is defined in (2).
\end{description}
Then (2) holds.
}

\medskip
\noindent
{\bf Proof.} See Appendix A.

\bigskip

%%%%%%%%%%%%%%%%%%%%%%%%%%%%%%%%%%%%%%%%%%%%%%%%%%%%%%%%%%%%%%%%%%%%
                  % R E M A R K      3.1
%%%%%%%%%%%%%%%%%%%%%%%%%%%%%%%%%%%%%%%%%%%%%%%%%%%%%%%%%%%%%%%%%%%%

\noindent {\bf Remark 3.1}

%%%%%%%%%%%%%%%%%%%%%%%%%%%%%%%%%%%%%%%%%%%%%%%%%%%%%%%%%%%%%%%%%%%%%%%%%%
 Before analyzing  the above results, let us  understand how the procedure works.
Consider the maximum likelihood recursive procedure in the
one-dimensional case
$$
\hat\theta_t=\hat \theta_{t-1}+I_t^{-1}(\hat \theta_{t-1})l_t(\hat
\theta_{t-1}),
$$
where $l_t(\theta) =\dot f_t^T(\theta, X_t \mid X_1^{t-1})/{f_ t(\theta,
X_t \mid X_1^{t-1})}$
  and $
I_t(\theta) $ is the conditional  Fisher information.

\vskip+0.1cm
\noindent
Denote  $\Dl_t=\hat\theta_t-\theta$
 and rewrite the above recursion as
$$
\Dl_t=\Dl_{t-1}+ I_t^{-1}(\theta+\Dl_{t-1}) l_t(\theta+\Dl_{t-1}).
$$
Then,
$$
E_\theta
\left\{\hat\theta_t-\hat\theta_{t-1}\mid{\cf}_{t-1}\right\}
=E_\theta \left\{\Dl_t-\Dl_{t-1}\mid{\cf}_{t-1}\right\} =
I_t^{-1}(\theta+\Dl_{t-1})b_t(\theta,\Dl_{t-1}),
$$
where
$$
b_t(\theta,u)= E_\theta \left\{l_t(\theta+u)\mid{\cf}_{t-1}\right\}.
$$
Under usual regularity conditions (see Sharia (2006a) Remark 3.2 for details),
$
b_t(\theta, 0)=0$ and $\frac{\partial}{\partial u} b_t(\theta, u)\mid_{u=0}=-i_t(\theta)<0,
$
implying that
\begin{equation}\label{mon}
u b_t(\theta,u) < 0
\end{equation}
for small values of $u\not=0$. Now, assuming that \eqref{mon} holds for all $u\not=0,$
suppose  that at time ~$t-1,$ ~
$\hat\theta_{t-1}<\theta,$ that is, $\Dl_{t-1}<0.$
 Then, by \eqref{mon},
$E_\theta
\left\{\hat\theta_t-\hat\theta_{t-1}\mid{\cf}_{t-1}\right\}>0.$
So, the next step $\hat\theta_t$ will be in the direction of
$\theta$.
 If  at time ~$t-1,$ ~ $\hat\theta_{t-1}>\theta,$
by the same reason, $E_\theta
\left\{\hat\theta_t-\hat\theta_{t-1}\mid{\cf}_{t-1}\right\}<0.$
So,  on average, at each step the
procedure moves towards $\theta$. However, the magnitude of the
jumps $\hat\theta_t-\hat\theta_{t-1}$ should decrease, for
otherwise, $\hat\theta_t$ may oscillate around $\theta$ without
approaching it. On the other hand,
care should be taken to  ensure that  the jumps do not decrease too rapidly  to avoid
failure of $\hat\theta_t$ to reach $\theta.$

\vskip+0.1cm
\noindent
These issues are addressed in Sharia (2006a)  and the conditions are introduced to ensure
global convergence of \eqref{rec2}, that is, convergence for any arbitrary starting value.
These conditions are flexible enough to allow for a quite wide choice of the normalising
sequence $\Gamma$ for any particular $\psi$.

\bigskip

%%%%%%%%%%%%%%%%%%%%%%%%%%%%%%%%%%%%%%%%%%%%%%%%%%%%%%%%%%%%%%%%%%%%
                  % R E M A R K      3.2
%%%%%%%%%%%%%%%%%%%%%%%%%%%%%%%%%%%%%%%%%%%%%%%%%%%%%%%%%%%%%%%%%%%%

\noindent {\bf Remark 3.2}

%%%%%%%%%%%%%%%%%%%%%%%%%%%%%%%%%%%%%%%%%%%%%%%%%%%%%%%%%%%%%%%%%%%%%%%%%%%%%%%%%%%
\noindent {\bf (i)}
As was mentioned above, strong consistency of the recursive estimator $\hat\theta_t$, that is
the convergence $\Dl_t=\hat\theta_t-\theta \to 0 $
($P^\theta$-a.s.) is established in Sharia (2006a).
Here we are interested in the asymptotic behaviour
of the recursive estimator given that it is consistent. Note
that although consistency is not formally required  in Lemma 3.1, it is
easy to see that if  $\hat\theta_t$ is not consistent, conditions (1)
and (2) will be satisfied for very special cases only.
Note also that given that $\Dl_t=\hat\theta_t-\theta \to 0$,
conditions (1) and (2) are  local in the sense  that they are determined
by local (w.r.t. the parameter) behaviour of the functions involved.

%%%%%%%%%%%%%%%%%%%%%%%%%%%%%%%%%%%%%%%%%%%%%%%%%%%%%%%%%%%%%%%%%%%%%%%%%%%%
\medskip
\noindent
 {\bf (ii)} Condition  (E) is an ergodicity
type assumption on the statistical model. If $\Gamma_t(\theta)=I_t(\theta)$
(the conditional Fisher information) and
$A_t(\theta)$ and $\eta(\theta)$ are non-random, then the model is called ergodic.
 Further discussion of this concept and related work
appears in  Basawa and  Scott
(1983), Hall and Heyde (1980) $\S$ 6.2, and
 Barndorff-Nielsen and  Sorensen (1994).

%%%%%%%%%%%%%%%%%%%%%%%%%%%%%%%%%%%%%%%%%%%%%%%%%%%%%%%%%%%%%%%%%%%%%%%%%%%%%%%%%
\medskip
\noindent {\bf (iii)}  Let us examine  condition   (2) in Lemma 3.1.
Given that $\Dl_t=\hat\theta_t-\theta \to 0$, if
 the functions $\p_t(\theta)$ and $\Gm_t(\theta)$  are continuous w.r.t.
$\theta$ (with certain uniformity w.r.t. $t$), we expect
${\cal E}_t(\theta) \to 0.$  Parts (b) and (c) in Proposition 3.1
give sufficient conditions for (2). If there exists a non-random sequence
$A_t(\theta),$ then obviously (L2) is less restrictive then
(LL2).  But unfortunately, (L2) can only be used for non-random
$A_t(\theta)$. In the case of random
$A_t(\theta)$,   when (LL2) may be used,   just the convergence
$E_\theta \left\{ \left({\cal E}_t(\theta)\right)^2
\mid{\cf}_{t-1}\right\} \to 0$ may not be enough since in many models the components
of $A_t(\theta)$ have the rate $\sqrt{t}$.
In such cases one may also use the result on the rate of convergence
of $\hat \theta_t$ presented in Sharia (2006b) (see examples 4.1 and 4.3 in the next section).

%%%%%%%%%%%%%%%%%%%%%%%%%%%%%%%%%%%%%%%%%%%%%%%%%%%%%%%%%%%%%%%%%%%%%%%%%%%%%%%
\medskip
\noindent {\bf (iv)}
 Condition (1)  gives an important clue for an optimal choice
 of the normalizing sequence $\Gamma_t(\theta)$.  To see this, let us assume that $\bf \psi \in \Psi^M$
so that $R_t(\theta,0)=0$  and have a look at (1) and (L1) in the  case of one dimensional parameter
$\theta \in \mathbb{R}.$ Now we can write
$$
\tr\Gm_t(\theta)\Dl_{t-1}+R_t(\theta,\Dl_{t-1})=
\left(\tr\Gm_t(\theta)+\frac {R_t(\theta,\Dl_{t-1}))-R_t(\theta,0)}{\Dl_{t-1}}\right)\Dl_{t-1}.
$$
In most applications,
the rate of $A_t$  is $\sqrt{t}$ and the best one can hope for is
that  $\sqrt{t} \Dl_t $ is stochastically bounded.
Therefore we must at least have the convergence
$\tr\Gm_t(\theta)+(R_t(\theta,\Dl_{t-1}))-R_t(\theta,0))/{\Dl_{t-1}}\to 0$. Given that
$\Dl_{t-1}\to 0$ we expect
$
\tr\Gm_t(\theta)\approx -{\partial}/{\partial u}~  R_t(\theta,u)\mid_{u=0}
$
for large $t$'s.
Also, since   $R_t(\theta,0)=
E_\theta \left\{\psi_t(\theta)
\mid{\cf}_{t-1}\right\}=0 $,  if  $\Gm_t(\theta)/\Gm_t(\theta+u)$ is smooth in $u=0$,
we can write that
$
{\partial}/{\partial u}~  R_t(\theta,u)\mid_{u=0} = \partial/\partial u ~E_\theta \left\{\psi_t(\theta+u)
\mid{\cf}_{t-1}\right\}\mid_{u=0}.
$
So, denoting
$$
 b_t(\theta,u)=E_\theta \left\{\psi_t(\theta+u)
\mid{\cf}_{t-1}\right\}
$$
we expect
%                             (DlGm)
\begin{equation}\label{DlGm}
\tr\Gm_t(\theta)\approx  -b'_t(\theta,0),
\end{equation}
where
$$
b'_t(\theta,0)=\frac{\partial}{\partial u} b_t(\theta, u)\mid_{u=0}.
$$
Using the similar arguments, for the multidimensional case,
we expect \eqref{DlGm} to hold for large $t$'s, where $b'_t(\theta,0)$ is
the total differential  of $b_t(\theta,u)$ in $u=0.$
Therefore,
%                             (Gm)
\begin{equation}\label{Gm}
\Gm_t(\theta)=-\sum_{s=1}^t b'_s(\theta,0)
\end{equation}
 is an obvious candidate for the normalizing sequence. If
$\p_t(\theta)$ is differentiable in $\theta$ and differentiation
of $ b_t(\theta,u)=E_{\theta} \{ \p_t(\theta+u)\mid {\cf}_{t-1}\}$
is allowed under the integral sign, then $b'_t(\theta,0)=E_{\theta} \{ \dot\p_t(\theta)\mid {\cf}_{t-1}\}.$
 This implies that,
for a given sequence of estimating functions  $\psi_t(\theta),$ another possible
 choice of the normalizing sequence is
%                             (Gm1)
\begin{equation}\label{Gm1}
\Gm_t(\theta)=-\sum_{s=1}^t E_{\theta} \{ \dot\p_s(\theta)\mid
{\cf}_{s-1}\},
\end{equation}
or any sequence with the increments
$$
\Dl \Gm_t=\Gm_t(\theta)-\Gm_{t-1}(\theta)= -E_{\theta} \{ \dot\p_t(\theta)\mid
{\cf}_{t-1}\}.
$$
Also, if the differentiation w.r.t. $\theta$ of
$$
%E_{\theta}\{ \p_t(\theta)\mid {\cf}_{t-1}\}
0=\int \p_t(\theta,z\mid X_1^{t-1}) f_t(\theta,z\mid X_1^{t-1})\mu
(dz)
$$
 is allowed under the integral sign, then by the product rule,
$$
0=\int \dot\p_t(\theta,z\mid X_1^{t-1})
f_t(\theta,z\mid X_1^{t-1})\mu (dz)+\int \p_t(\theta,z\mid
X_1^{t-1}) \dot f_t(\theta,z\mid X_1^{t-1})\mu (dz).$$ So,
\begin{eqnarray}\label{long}
 E_{\theta}\{ \dot\p_t(\theta)\mid {\cf}_{t-1}\}
& & =\int \dot\p_t(\theta,z\mid X_1^{t-1})
f_t(\theta,z\mid X_1^{t-1})\mu (dz) \nonumber \\
& &=-\int \p_t(\theta,z\mid X_1^{t-1})
\dot f_t(\theta,z\mid X_1^{t-1})\mu (dz) \nonumber \\
& &=-\int \p_t(\theta,z\mid X_1^{t-1})
 l^{T}_t(\theta,z\mid X_1^{t-1})f_t(\theta,z\mid X_1^{t-1})\mu (dz) \\
&&=-E_{\theta}\{ \p_t(\theta)l^{T}_t(\theta)\mid
{\cf}_{t-1}\},\nonumber
\end{eqnarray}
where, as before, $l_t(\theta)=\dot f_t^T(\theta, X_t|X_1^{t-1})/f_t(\theta, X_t|X_1^{t-1}).$
Therefore, denoting
$$
\gm_t^\p(\theta)=E_{\theta}\{ \p_t(\theta)l^{T}_t(\theta)\mid
{\cf}_{t-1}\},
$$
another possible choice of the normalizing sequence is
%                             (Gm2)
\begin{equation}\label{Gm2}
\Gm_t(\theta)=\sum_{s=1}^t\gm_s^\p(\theta),
\end{equation}
or any sequence with the increments
$$
\Dl \Gm_t=\Gm_t(\theta)-\Gm_{t-1}(\theta)= \gm_t^\p(\theta).
$$
Since typically, for each $\theta,$ the process
$$
M_t^\theta=\sum_{s=1}^t\p_s(\theta)
$$
is a $P^\theta$ -- martingale, \eqref{Gm2} can be rewritten as
$$
\Gm_t(\theta)=\langle M^\theta, U^\theta \rangle_t
$$
where $U_t^\theta=\sum_{s=1}^t l_s(\theta)$  is the  score martingale.

%%%%%%%%%%%%%%%%%%%%%%%%%%%%%%%%%%%%%%%%%%%%%%%%%%%%%%
\medskip
\noindent {\bf (v)}
Part (iv) above highlights a very important point. Suppose  we wish to construct
a recursive estimator with a given sequence $\psi$  of estimating functions.
In order to achieve consistency, we are quite flexible in choice of the  normalizing
sequence $\Gamma$;
the recursive procedure will converge even when $\Gamma$ sequence is not related to
$\psi$ (see Sharia (2006a)). (Of course, the rate of the normalizing sequence  still has to be ``right''
but is mostly determined by the model.)
If we want to obtain a recursive estimator which is also
 asymptotically  linear, then the  normalizing sequence $\Gamma$ has to be
\eqref{Gm}  (or \eqref{Gm1},  \eqref{Gm2}, or a sequence asymptotically equivalent to
 \eqref{Gm}).

%%%%%%%%%%%%%%%%%%%%%%%%%%%%%%%%%%%%%%%%%%%%%%%%%%%%%%%%%%%%%%%%%%%%%%%%%%%%%%%%%%%%%%%%%%%%%%%%%%%
\medskip
\noindent {\bf (vi)} Let us consider a likelihood case,
 that is   $\psi_t(\theta)=l_t(\theta).$
Since $\gm_t^\p(\theta)=i_t(\theta),$ the process  \eqref{Gm2} in this case is the conditional Fisher
information $I_t(\theta)=\sum_{s=1}^ti_s(\theta).$ So, the corresponding recursive procedure is
%%%%%%%%%%%%%%%%%%%%%%%%%%%%%%%%%%%%%%%%%%%%%%%%%%
%                   \eqref{recmle}
%%%%%%%%%%%%%%%%%%%%%%%%%%%%%%%%%%%%%%%%%%%%%%%%%
\begin{equation} \label{recmle}
\hat\theta_t=\hat\theta_{t-1}+
{I_t^{-1}(\hat\theta_{t-1})}l_t(\hat\theta_{t-1}),
~~~~~~~~~t\ge 1,
\end{equation}
Also, given that the model possesses certain  ergodicity properties,
asymptotic linearity of \eqref{recmle} implies asymptotic
efficiency.   In particular, in the case
of i.i.d. observations, it follows that the above recursive procedure is
asymptotically
normal with parameters $(0, \; i^{-1}(\theta) )$ (see Corollary 4.1 in Section 4).

%%%%%%%%%%%%%%%%%%%%%%%%%%%%%%%%%%%%%%%%%%%%%%%%%%%%%%%%%%%%%%%%%%%%%%%%%%%%%%%%%%%%%%%%%%%%%%%%%%%
\medskip
\noindent {\bf (vii)} Normalizing sequences suggested   in (iv)
have been derived from the asymptotic considerations. In practice however, behaviour of $\Gamma$
sequence  for the first several steps might also  be important. This can happen when the number of observations
is small or even moderately large.
According to (iv), to achieve asymptotic linearity, one has to choose a normalizing sequence $\Gamma$
with the property that
$$
\tr\Gm_t(\theta)\approx -b'_t(\theta,0)
$$
for large $t$'s. So, we can consider any sequence of the form $C+c_t\Gamma_t$, where $\Gamma_t$ is one of the sequences introduced
above (by \eqref{Gm}, \eqref{Gm1}, or \eqref{Gm2}),
 $c_t$ is a sequence of non-negative r.v.'s such that $c_t =1$  eventually
  and $C$ is a suitably chosen constant.
 In practice,  $c_t$ and $C$  can be treated as tuning constants
  to control behaviour
of the procedure for the first several steps (see Sharia (2006a), Remark 4.4).
Under certain assumptions, at each step, the recursive procedure \eqref{rec2},
(on average) moves towards the direction of the unknown parameter (see Remark 3.1 or  Sharia (2006a), Remark 3.2
for details).
Nevertheless, if the values of the normalizing sequence are too small  for the first several steps, then the
procedure will oscillate excessively around the true value of the parameter. On the other hand, too
large  values of the normalizing sequence will result in slower  convergence of the procedure. A good
balance can be achieved by using the tuning constants.
The detailed discussion of these and related topics will appear elsewhere, but as a rough guide,
the graph of $\hat\theta_t$ against $t$ should ideally have a shape of
those in Figure 1 in Sharia (2006a) (that is, a reasonable  oscillation at the beginning  of
the procedure before settling down at a particular level).

%\vskip+0.5cm

                % 4  SPECIAL MODELS   AND EXAMPLES

\section{SPECIAL MODELS
                   AND EXAMPLES}
\setcounter{equation}{0}

                      %  E X A M P L E 1
\noindent {\large \bf 4.1.} {\bf  The i.i.d. scheme.} Consider  the
classical scheme of
  i.i.d. observations $X_1,X_2,\ldots ,$ with a common probability
density/mass  function $f(\theta,x), \;\; \theta \in
{\mathbb{R}}^m.$ Suppose that $\p(\theta, x)$ is an estimating
function with
$$
E_\theta(\p(\theta,X_1))=\int \p(\theta,z)f(\theta,z)\mu (dz)=0.
$$
Let us define the recursive estimator $\hat \theta_t$  by
                             % mleg
\begin{equation}\label{mleg}
\hat \theta_t = \hat \theta_{t-1}  + \frac 1 t\gm^{-1} (\hat
\theta_{t-1})
  \p(\hat \theta_{t-1},X_t),\qquad t\ge 1,
\end{equation}
where $\hat\theta_0\in {\mathbb{R}}^m$ is any initial value.
According to Remark 3.2 (iv) and the condition (V) below,
an optimal choice of $\gm(\theta)$ would be either
$$
\gm(\theta)=E_\theta(\dot\p(\theta,X_1))
 %\int\dot\p(\theta,z) f(\theta,z)\mu(\,dz),
$$
or
$$
\gm(\theta)=
% \int\p(\theta,z) l^T(\theta,z)f(\theta,z)\mu(\,dz),
 E_\theta(\p(\theta,X_1) l^T(\theta,X_1)) ~~~~~  \mbox{where} ~~~~~
l(\theta,x)=\frac{\dot f^T(\theta,x)}{f(\theta,x)},
$$
 or any  non-random invertible  matrix function that satisfies conditions
 listed below.

Suppose that
$$
   j_\p(\theta)=\int
\p(\theta,z)\p^T(\theta,z) f(\theta,z)\mu(\,dz)< \infty
$$
and consider the following conditions.
\begin{description}
\item[(I)] For any $0<\ve<1,$
$$
\sup_{\ve\le \|u\|\le \frac 1{\ve}}
u^T \;\;  \gm^{-1}(\theta+u)\int \p(\theta+u,x)f(\theta,x)\mu(\,dx)  <  0.
$$
\item[(II)] $\;$For each $ u \in {{\mathbb{R}}^m},$
$$
\int \left \|  \gm^{-1}(\theta+u) \p(\theta+u,x)\right \|^2
f(\theta,x)\mu(\,dx) \leq  K_\theta
(1+\|u\|^2)
$$
for some constant $K_\theta.$
\item[(III)] $\;$
$\gm(\theta)$ is continuous in $\theta.$
\item[(IV)]
$$
\lim_{u\to 0} \int {\|\p(\theta+u,x)-\p(\theta,x)\|}^2
f(\theta,x)\mu (\, dx) = 0.
$$
\item[(V)]
$$
\int \p(\theta+u,x)f(\theta,x)\mu(\,dx)=-\gm(\theta+u) u +\alpha^\theta(u),
$$
where $\alpha^\theta(u)=o(\|u\|^{1+\ve})$ as  $u\to 0$ for some $\ve>0$.
\end{description}

%%%%%%%%%%%%%%%%%%%%%%%%%%%%%%%%%%%%%%%%%%%%%%%%%%%%%%%%%%%%%%%%
%                             C o r  4.1
%%%%%%%%%%%%%%%%%%%%%%%%%%%%%%%%%%%%%%%%%%%%%%%%%%%%%%%%%%%%%%%%%

                % C O R O L L A R Y   4.1
\begin{cor}
Suppose that  for any $\theta \in {\mathbb{R}}^m$ conditions {\bf
(I)} - {\bf (V)} are satisfied. Then  the   estimator $\hat
\theta_t$   is strongly consistent and ~ $ t^\dl(\hat
\theta_t-\theta) \to 0  ~(P^\theta$-a.s.) for any $0<\dl<1/2$ and
any initial value $\hat\theta_0$. Furthermore, $\hat \theta_t$ is
asymptotically normal with parameters $(0, \;\; \gm^{-1}(\theta)
j(\theta,0) \gm^{-1}(\theta) )$, that is,
$$
{\cal L}\left( t^{1/2}(\hat  \theta_t -\theta) \mid P^\theta
 \right) {\stackrel{w}{\to}} {\cal N} \left(0,\;\; \gm^{-1}(\theta)
 j_\p(\theta) \gm^{-1}(\theta)\right).
$$
In particular, in the  case of  the maximum likelihood type recursive procedure
with $\psi(\theta,x)=\dot f^T(\theta,x)/f(\theta,z)$ and
$\gamma(\theta)=i(\theta)=j_l(\theta)$, the estimator $\hat
\theta_t$ is asymptotically efficient  (i.e., asymptotically
normal with parameters $(0, \;\; i^{-1}(\theta) )$).
\end{cor}

%%%%%%%%%%%%%%%%%%%%%%%%%%%%%%%%%%%%%%%%%%%%%%%%%%5555555555555555
\noindent {\bf Proof} See Appendix A.

%%%%%%%%%%%%%%%%%%%%%%%%%%%%%%%%%%%%%%%%%%%%%%%%%%5555555555555555

\bigskip

Similar  results (for i.i.d. schemes) were obtained
 by  Khas'minskii and Nevelson  (1972) (when $\p(\theta,x)=l(\theta,x)$ and
 $\gm(\theta)=i(\theta)$,
Ch.8, $\S$4) and Fabian (1978).

%%%%%%%%%%%%%%%%%%%%%%%%%%%%%%%%%%%%%%%%%%%%%%%%%%%%%%%%%%%%%
%       Linear recursion
%%%%%%%%%%%%%%%%%%%%%%%%%%%%%%%%%%%%%%%%%%%%%%%%%%%%%%%%%%%%%
                     %  E X A M P L E 2
\noindent {\large \bf 4.2.} {\bf Linear  procedures.} Consider the
recursive procedure
%          LinRecg
\begin{equation}\label{LinRecg}
\hat\theta_t=\hat\theta_{t-1}+\Gamma_t^{-1}\left( h_t-
\gamma_t\hat\theta_{n-1}\right), \;\;\; t\ge 1,
\end{equation}
where the $\Gamma_t$ and  $\gamma_t$ are predictable matrix processes,
$h_t$ is an adapted  process (i.e., $h_t$ is
$\mathcal{F}_t$-measurable for $t\ge 1$)
 and all three are independent of
$\theta.$ The following result  gives a sets of sufficient
conditions for the asymptotic linearity  of the estimator defined by \eqref{LinRecg}
 in the case when
the linear $\psi_t(\theta)=h_t-\gamma_t\theta$ is a
martingale-difference, i.e.,
 $E_\theta\left\{h_t\mid {\cal{F}}_{t-1}\right\}
=\gamma_t\theta,$ for $t \ge 1. $

%%%%%%%%%%%%%%%%%%%%%%%%%%%%%%%%%%%%%%%%%%%%%%%%%%%%%%%%%%%%%%%%
%                             C o r  4.2
%%%%%%%%%%%%%%%%%%%%%%%%%%%%%%%%%%%%%%%%%%%%%%%%%%%%%%%%%%%%%%%%%

                % C O R O L L A R Y   4.2
\begin{cor}
 Suppose that  $\Gamma_t \to \infty $  and
                             % mleg
\begin{equation}\label{L1Lin}
\Gamma_t^{-1/2}\sum_{s=1}^t(\tr\Gm_s-\gamma_s)\Dl_{s-1} \to 0
\end{equation}
in probability $P^\theta$, where $\Dl_{s-1}=\hat\theta_{s-1}-\theta.$
Then the recursive estimator defined by \eqref{LinRecg} is asymptotically
linear with
\begin{equation}
\Gm^{1/2}_t(\hat \theta_t -\theta) = \Gm^{-1/2}_t\sum_{s=1}^t
\p_s(\theta)+o_{P^\theta}(1),
\end{equation}
where  $o_{P^\theta}(1)\to 0$  in probability  $P_\theta.$
\end{cor}

%%%%%%%%%%%%%%%%%%%%%%%%%%%%%%%%%%%%%%%%%%%%%%%%%%5555555555555555
\noindent {\bf Proof}
Let us check the conditions of Lemma 3.1 for $A_t(\theta)=\Gamma_t^{1/2}.$
Condition (E) trivially holds. Then, since $\psi_t(\theta)=h_t-\gamma_t\theta$ and
$$
b_t(\theta,u)=E_\theta\left\{(\psi_t(\theta+u))\mid
{\cal{F}}_{t-1}\right\}=E_\theta\left\{(h_t-\gamma_t(\theta+u))\mid
{\cal{F}}_{t-1}\right\}=-\gamma_tu,
$$
we have
$$
R_t(\theta,u)=
                    \Gm_t(\theta)
            \Gm_t^{-1}(\theta+u)b_t(\theta,u)
            =-\gamma_tu.
            $$
Therefore, (1) is equivalent to \eqref{L1Lin}.
Then, it is easy to see that  for ${\cal E}_s(\theta)$  defined in (2) we have
$$
{\cal E}_s(\theta)=
\p_s(\theta+\Dl_{s-1})-b_s(\theta,\Dl_{s-1})
-\p_s(\theta)=0
$$
implying that (2) holds which completes the proof.
 $\diamondsuit $
\bigskip

%%%%%%%%%%%%%%%%%%%%%%%%%%%%%%%%%%%%%%%%%%%%%%%%%%%%%%%%%%%%%%%%%%%%%%%%%%
                  % R E M A R K 4.1
%%%%%%%%%%%%%%%%%%%%%%%%%%%%%%%%%%%%%%%%%%%%%%%%%%%%%%%%%%%%%%%%%%%%%%%%%%
\begin{rem}
{\rm
Condition \eqref{L1Lin} trivially holds if
$\Delta\Gamma_t=\gamma_t,$ that is  $\Gamma_t=\sum_{s=1}^t \gamma_s.$
In this case, the solution of \eqref{LinRecg} is
%         Solg
  \begin{equation}\label{Solg}
\hat\theta_t=\Gamma_t^{-1}\left(\hat\theta_0+\sum_{s=1}^t
h_s(X_s)\right).
\end{equation}
This can be easily seen by inspecting the difference
$\hat\theta_t-\hat\theta_{t-1}$ for the sequence \eqref{Solg} (exactly in the same way as
in \eqref{solut}), to
check that \eqref{LinRecg} holds.  Also, since  \eqref{Solg} can obviously be rewritten as
$$
\hat\theta_t=\Gamma_t^{-1}\hat\theta_0+\Gamma_t^{-1}\sum_{s=1}^t
\left(h_s(X_s)- \gamma_s\theta\right) +\theta,
$$
it follows that in this case,
$\Gamma_t \to \infty$ is indeed
 an obvious  necessary and sufficient condition for $\hat\theta_t$ to be
 asymptotically linear (for arbitrary  starting value $\hat\theta_0$).
 }
\end{rem}
\vskip+0.5cm

\vskip+0.3cm

                      %  E X A M P L E 3

\noindent {\large \bf 4.3.}
 {\bf  Exponential family of  Markov
processes} Consider a conditional exponential family of  Markov
processes in the sense of Feigin (1981) (see also Barndorf-Nielson
(1988)). This is a time homogeneous Markov chain
with the one-step transition density
$$
f(y; \theta,x)=h(x,y)\exp\left(
\theta^Tm(y,x)-\beta(\theta;x)
\right),
$$
where $m(y,x)$ is a $m$-dimensional vector and $\beta(\theta;x)$
is one dimensional. Then in our notation
$f_t(\theta)=f(X_{t};\theta ,X_{t-1})$ and
$$
l_t(\theta)=\left(\frac{d}{d\theta}\log f_t(\theta)\right)^T=
m(X_t,X_{t-1})-\dot\beta^T(\theta;X_{t-1}).
$$
 It follows from standard
exponential family theory (see, e.g., Feigin (1981)) that
$l_t(\theta)$
is a martingale-difference and
 the  conditional Fisher information
is
$$
I_t(\theta)=\sum_{s=1}^t\ddot\beta(\theta;X_{s-1}).
$$
A maximum likelihood  type recursive procedure  can be defined as
$$
\hat\theta_t=\hat\theta_{t-1}+\left(\sum_{s=1}^t\ddot\beta(\hat\theta_{t-1};X_{s-1})
\right)^{-1}\left(m(X_t,X_{t-1})-\dot\beta^T(\hat\theta_{t-1};X_{t-1})
\right), ~~~ t\ge  1.
$$

Now suppose that $\theta$ is one dimensional and the process belongs to the
conditionally additive exponential family, that is,
$$
f(y; \theta,x)=h(x,y)\exp\left(
\theta m(y,x)-\beta(\theta;x)
\right),
$$
with
\begin{eqnarray}\label{Add}
\beta(\theta;x)=\gamma(\theta) h(x)
\end{eqnarray}
where $h(\cdot)\ge 0$  and $\ddot \gamma(\cdot) \ge 0$
(see Feigin (1981)).
Then,
$$
I_t(\theta)=\ddot \gamma(\theta)H_t
~~~ \mbox{where} ~~~
H_t=\sum_{s=1}^t h(X_{s-1}).
$$
Assuming that $\ddot \gamma(\theta)\not= 0,$
the likelihood recursive procedure is
\begin{eqnarray}\label{RecAdd}
\hat\theta_t=\hat\theta_{t-1}+\frac1{\ddot\gamma(\hat\theta_{t-1})H_t}
\left(m(X_t,X_{t-1})-\dot\gamma(\hat\theta_{t-1})h(X_{t-1})
\right).
\end{eqnarray}

%%%%%%%%%%%%%%%%%%%%%%%%%%%%%%%%%%%%%%%%%%%%%%%%%%%%%%%%%%%%%%%%%%%%%%%%%%
                  % R E M A R K 4.2
%%%%%%%%%%%%%%%%%%%%%%%%%%%%%%%%%%%%%%%%%%%%%%%%%%%%%%%%%%%%%%%%%%%%%%%%%%
\begin{rem}
{\rm
 Consistency and rate of convergence of the estimator derived by \eqref{RecAdd} is
 studied  In Sharia (2006b).
To ensure that  \eqref{RecAdd} has the same asymptotic properties as the maximum likelihood
estimator, one has to impose certain restrictions  on the $\gamma(\theta)$ and $H_t$.
In Corollary A1 in Appendix A,  the  conditions of Section 3 written in terms of   this model
are presented. These conditions  will be satisfied if there is a certain balance between requirements of
  smoothness on $\gamma(\cdot)$,  the rate at which $H_t \to \infty$, and  ergodicity
of the model. For instance, suppose that  the model is ergodic, that is, there exists a non-random sequence
$\tilde H_t$ such that $H_t/\tilde H_t \to \eta < \infty$ weakly. Then
$$
\frac 1{I_t^{1/2}(\theta)}\sum_{s=1}^t {\cal E}_s(\theta) \to  0,
$$
will hold if the process
$$
\frac 1{I_t(\theta)}\sum_{s=1}^t E_\theta \left\{
{\cal E}_s^2(\theta) \mid{\cf}_{s-1}\right\}
=\frac 1{I_t(\theta)}\sum_{s=1}^t\tr I_s(\theta)
\left(\frac{\ddot\gamma(\theta+\Dl_{s-1})-\ddot\gamma(\theta)}{\ddot\gamma(\theta+\Dl_{s-1})}\right)^2,
$$
 converges to zero (criterion based on the Lenglart-Rebolledo inequality, see (L2) and  formula (A5) in Appendix
 A). So, assuming that
 the estimator is consistent (that is $\Dl_t\to 0$),  by the Toeplits lemma, the above will be guaranteed
 by the continuity of $\ddot \gamma_t(\cdot)$. On the other hand, if
 the model is non-ergodic, then one may need to impose  smoothness of higher order on $\gamma(\cdot)$
   function (see condition (iii) below) and restrictions on the growth
 of the sequence $H_t$ (see condition (i) below). The following result gives one possible
 set of sufficient conditions for the recursive estimator to be consistent and to have the same
 asymptotic properties  as the maximum likelihood estimator.
}
 \end{rem}

\bigskip
\noindent
{\bf Proposition  4.3}
{\it
Suppose that $H_t\to \infty$ and
\begin{description}
\item[(i)]
$$
  \frac{h(X_{t})}{H_{t}}\to 0;
 $$
\item[(ii)] there exists a constant $B$ such that
$$
\frac{1+\dot\gamma^2(u)}{\ddot\gamma^2(u)} \le B(1+u^2)
$$
for each $ u\in {\mathbb{R}}$.
\item[(iii)] The function $\ddot\gamma(\cdot)$ is locally
Lipschitz , that is, for any $\theta$ there exists
a constant $K_\theta$ and $0<\ve_\theta \le 1/2$ such that
$$
|\ddot\gamma(\theta+u)-\ddot\gamma(\theta)|\le K_\theta |u|^{\ve_\theta}
$$
for small  $u$'s.
\end{description}
Then
 $\hat\theta_t$  defined by \eqref{RecAdd} is  strongly consistent (i.e.,
$\hat \theta_t \to \theta \;\; P^\theta$-a.s.) for any  initial
value $\hat \theta_0$. Furthermore,   $H_t^\delta(\hat \theta_t - \theta )\to 0\;\;
P^\theta$-a.s. for any $\delta\in ]0,1/2[$,
and $\hat\theta_t$  is asymptotically linear  with
\begin{equation}\label{aslinexp}
H^{1/2}_t(\hat \theta_t -\theta) = H^{-1/2}_t \sum_{s=1}^t
\left(m(X_s,X_{s-1})-\dot\gamma(\theta)h(X_{s-1})\right)
+o_{P^\theta}(1),
\end{equation}
where  $o_{P^\theta}(1)\to 0$  in probability  $P_\theta.$

}

\bigskip

                      %  E X A M P L E 4

\noindent {\large \bf 4.4.}  {\bf  AR(m) process}
 Consider an  AR(m) process
$$
X_i=\theta_{1} X_{i-1}+\dots+\theta_{m} X_{i-m}+\xi_i =\theta^T
X_{i-m}^{i-1}+\xi_i,
$$
where $X_{i-m}^{i-1}=(X_{i-1},\dots,X_{i-m})^T,$
$\theta=(\theta_1,\dots, \theta_m)^T$
 and
${\xi}_i$ is a sequence of i.i.d. random variables.

In Sharia (2006a) we discuss convergence of the recursive estimators of the form
                          %(Arg)
\begin{equation}\label{Arg}
\hat \theta_t=\hat \theta_{t-1}+\Gm_t^{-1}(\hat\theta_{t-1})
 \p_t(X_t-\hat \theta_{t-1}^TX_{t-m}^{t-1}),
\end{equation}
where $\p_t(z)$ and $\Gm_t^{-1}(z)$ ($z\in {\mathbb{R}}^m$) are
respectively  suitably chosen vector and matrix processes.
If  the probability density function of   ${\xi}_t$ w.r.t. Lebesgue's measure is
$g(x)$ then the conditional probability density function of $X_t$ given
values of past observations  of $X_{t-m}^{t-1}=(X_{t-1},\dots , X_{t-m})$ is
obviously
$$
f_ t(\theta, x_t \mid x_{t-m}^{t-1})=g(x_t-\theta^Tx_{t-m}^{t-1}),
$$
and so,
$$
l_t(\theta)=\frac{\dot f_ t^T(\theta, X_t \mid X_{t-m}^{t-1})}{f_ t(\theta, X_t \mid X_{t-m}^{t-1})}=
-\frac{ g' (X_t-\theta^T X_{t-m}^{t-1})}{g(X_t-\theta^T X_{t-m}^{t-1})}X_{t-m}^{t-1}.
$$
It follows from the results of Section 3 (see Remark 3.2 (vi)) that an optimal
choice of the normalizing sequence is the conditional Fisher information $I_t(\theta)$,
(or any sequence  with the increments equal to $\Dl I_t(\theta)$).
It is easy to see that in this case,
$$
 I_t(\theta)=I_t={i^g}
\sum_{s=1}^t X_{t-m}^{t-1}(X_{t-m}^{t-1})^T
$$
where
$$
i^g=\int  \left(\frac{\dot g'(z)}{g(z)}\right)^2
 g(z) \,dz.
$$
Since in this case the conditional Fisher information can also be found recursively,
a likelihood recursive procedure is
\begin{eqnarray}\label{ArLsq}
&&\hat \theta_t=\hat \theta_{t-1}-I_t^{-1}
\frac{ g' (X_t-\hat \theta_{t-1}X_{t-1})}{g(X_t-\hat \theta_{t-1}X_{t-1})}X_{t-m}^{t-1}
\\
&&I_t=I_{t-1}+i^gX_{t-m}^{t-1}(X_{t-m}^{t-1})^T, \nonumber
\end{eqnarray}
 for $t\ge 1$  and  an  arbitrary starting point $\hat\theta_0$.
The strong consistency of the estimators \eqref{Arg} and, in particular,
that of \eqref{ArLsq} is studied in Sharia (2006a).

The class of estimators \eqref{Arg} includes
 recursive versions of robust modifications of the
least squares method. These are recursive estimators defined by
\begin{equation}\label{Camb}
\hat \theta_t=\hat \theta_{t-1}+\Gamma_t\gm(X_{t-m}^{t-1})
 \phi(X_t-\hat \theta_{t-1}^TX_{t-m}^{t-1}),
\end{equation}
where
$\phi$ is a bounded  scalar function and $\gm(u)$ is a vector
function of the form $u h(u)$ for some non-negative function $h$
of $u.$

Since \eqref{Camb} is of the form \eqref{rec2} with
$\psi_t(\theta)=\gm(X_{t-m}^{t-1})\phi\left(X_t-\theta^TX_{t-m}^{t-1}\right),$
assuming  that $\phi(\cdot)$ is  differentiable (almost everywhere w.r.t. Lebesgue's measure)
we obtain
\begin{eqnarray}
E_\theta\left\{
\dot \psi_t(\theta)\mid {\cal{F}}_{s-1}\right\}
& =&-\gm(X_{t-m}^{t-1})(X_{t-m}^{t-1})^T
E_\theta\left\{
\phi'\left(X_t-\theta^TX_{t-m}^{t-1}\right)\mid {\cal{F}}_{s-1}
\right\} \nonumber\\
 &=&-\gm(X_{t-m}^{t-1})(X_{t-m}^{t-1})^T\int \phi'\left(x-\theta^TX_{t-m}^{t-1}\right)g(x-\theta^T X_{t-m}^{t-1})dx, \nonumber
\nonumber\\
 &=&-\gm(X_{t-m}^{t-1})(X_{t-m}^{t-1})^T\int \phi'(x)g(x)dx. \nonumber
\end{eqnarray}
 So, according to Lemma 3.1 (see Remark  3.2 (iv) formula \eqref{Gm1}),
an optimal normalizing sequence $\Gamma_t$ for \eqref{Camb}
is
%%%%%%%%%%%%%%%%%%%%%%%%%%%%%%%%%%%%%%%%%%%%%%%%%%%%%%%%%%%%
\begin{equation}\label{GG}
\Gamma_t(\theta)=C_g
\sum_{s=1}^t\gm(X_{s-m}^{s-1}) {X_{s-m}^{s-1}}^T
\end{equation}
where
$$
C_g=\int\phi'(x)g(x)dx
$$
or a sequence with the increments equal to
$C_g\gm(X_{s-m}^{s-1}) {X_{s-m}^{s-1}}^T.$

Consider for instance a  recursive M-estimator
of the parameter of  an AR(1) process  defined as
                          %(RecRob1)
\begin{equation}\label{RecRob}
\hat \theta_t=\hat \theta_{t-1}+\frac 1{\Gamma_t} s_x
\phi_c\left(\frac{X_{t-1}}{s_x}\right)
 s_r\phi_c\left(\frac{X_t-\hat \theta_{t-1}X_{t-1}}{s_r}\right)
\end{equation}
where $s_x$ and  $s_r$  are scale estimates and $\phi_c$ is the Huber function,
\begin{equation}
\phi_c(x)=\left\{\begin{array}{ll}
  x,   &\mbox{if $|x|\le c$}\\\nonumber
 c ~ \mbox{sign}(x)            &\mbox{if $|x|>c$}
\end{array}
\right.
\end{equation}
and $c>0$ is a tuning  constant. This is a
recursive version of a robust generalized M-estimator
of the parameter of  an AR(1) process  proposed by see Denby and Martin (1979).

Another example is
                          %(RecRob2)
\begin{equation}\label{RecRob2}
\zeta_t=\zeta_{t-1}+\frac 1{\Gamma_t^\zeta} s_x
\phi_{\alpha, \beta}\left(\frac{X_{t-1}}{s_x}\right)
 s_r\phi_{\alpha, \beta}
 \left(\frac{X_t-\zeta_{t-1}X_{t-1}}{s_r}\right).
\end{equation}
 where  $\phi_{\alpha, \beta}$
is  Hampel's two-part redescending  function
\begin{equation}
\phi_{\alpha,\beta}(x)=\left\{\begin{array}{llll}
  x,                       &\mbox{if $|x|\le \alpha$}\\
 \alpha(\beta-x)/(\beta-\alpha),            &\mbox{if $\alpha <x \le \beta$}\\
 -\alpha(\beta+x)/(\beta-\alpha),            &\mbox{if $-\beta \le x < -\alpha$}\\
   0,                       &\mbox{if $|x| \ge \beta$},
\end{array}
\right.
\end{equation}
with tuning constants $0 <\alpha <\beta.$

For the procedure \eqref{RecRob},
$$
C_g=\int\phi'(x)g(x)dx=\int s_r
\left(\frac{d}{dx}\phi_c\left(\frac{x}{s_r}\right)\right)g(x)dx=\int
\phi'_c\left(\frac{x}{s_r}\right)g(x)dx,
$$
and so
%%%%%%%%%%%%%%%%%%%%%%%%%%%%%%%%%%%%%%%%%%%%%%%%%%%%%%%%%%%%%%%%%%%%%%%%%%%%
%
\begin{eqnarray}\label{Cg1}
C_g=\int_{-cs_r}^{cs_r} g(x)dx
\end{eqnarray}
Similarly, for \eqref{RecRob2},
$$
C_g=\int
\phi'_{\alpha,\beta}\left(\frac{x}{s_r}\right)g(x)dx
$$
%%%%%%%%%%%%%%%%%%%%%%%%%%%%%%%%%%%%%%%%%%%%%%%%%%%%%%%%%%%%%%%%%%%%%%%%%%%%
%
\begin{equation}\label{Cg2}
=\int_{-\alpha s_r}^{\alpha s_r} g(x)dx
-\frac{\alpha}{\beta-\alpha}\left(\int_{-\beta s_r}^{-\alpha s_r} g(x)dx+\int_{\alpha s_r}^{\beta s_r}
g(x)dx\right)
\end{equation}

 Below we present a brief simulation study.
The time series were generated from the additive effect outliers
(AO) model:
\begin{eqnarray}
Y_t=\theta Y_{t-1}+w_t \nonumber \\
X_t=Y_t+v_t, \nonumber
\end{eqnarray}
where innovations  $w_t$ are i.i.d. Gaussian $N(0,1).$ The
variables  $v_t$ are also i.i.d. with  distribution $
(1-\ve)\dl_0+\ve N(0,\sg^2),$  where $\dl_0$ is the distribution
that assigns probability $1$ to the origin.  Therefore, with
probability $1-\ve$ the $AR(1)$ process $Y_t$ is observed, and
with probability $\ve$ the observation is the $AR(1)$ process
$Y_t$ plus the error with Gaussian distribution $N(0,\sg^2)$. In
this simulation, $\theta=0.6 $, $\ve=0.05$ and $\sg^2=9$. The
figures below show the performances of the estimator
$\hat\theta_t$ defined by \eqref{RecRob}, the estimator $\zeta_t$
defined by \eqref{RecRob2} and  the least squares estimator
$\hat\theta_t^{ls}$ (which is equivalent to the recursive
procedure defined by \eqref{ArLsq} with $\dot g(x)/g(x)=-x$). The
estimators are computed for the series of length $200 $, with the
additional $30$ observations
 at the beginning on which initial estimates are based; as an
 estimates for $s_x$ and $s_r$ we take the median of the absolute
 values of the data and  residuals respectively,  divided by
 0.6745. The p.d.f. $g(x)$ in \eqref{RecRob} and  \eqref{RecRob2}
is replaced by the p.d.f. of $N(0,s_r^2)$
 and the values of the tuning constants are
$c=1.8,$ $\alpha=1.8$ and $\beta=4$.
 Figure 1 shows single realizations   and the mean squared errors over 300 replications
 of the  estimators
$\hat\theta_t^{ls},$ $\hat\theta_t$ and  $\zeta_t$ for $t=5,\dots, 200$.

%\newpage

\begin{figure}
\begin{center}
\resizebox{0.9\textwidth}{!}{\includegraphics{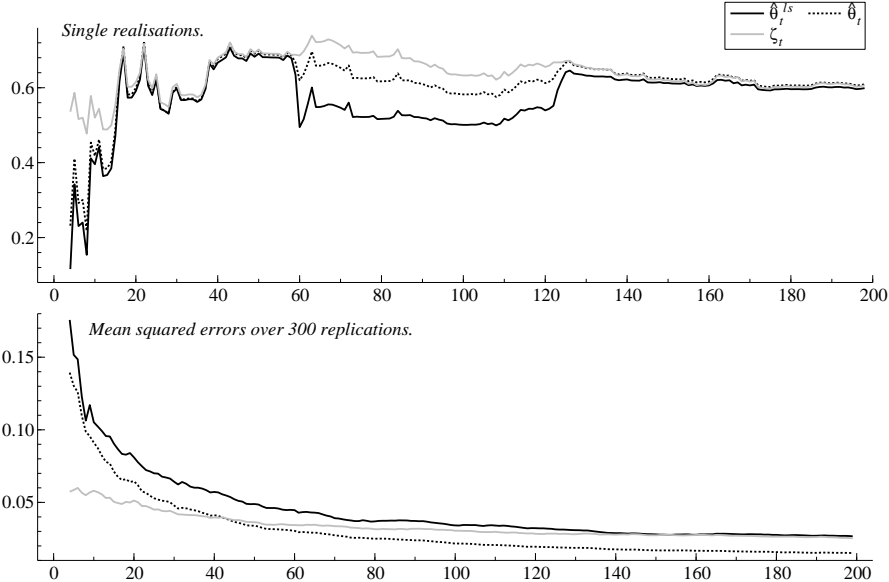}}
\end{center}
\caption{{\small  Single realizations   and the mean squared errors over 300 replications,
for $t=5,\dots, 200$.}}
\end{figure}
%
%\begin{figure}
%\begin{center}
%\resizebox{0.9\textwidth}{!}{\includegraphics{Ms}}
%\end{center}
%\end{figure}

%\end{figure}
%\begin{figure}
%\begin{center}
%\resizebox{0.9\textwidth}{!}{\includegraphics{MsGrey1}}
%\end{center}
%\caption{\small }

Further  simulation study is required to study performances of
these procedures.  As this brief simulation suggests, both
$\hat\theta_t$  and $\zeta_t$ outperform  $\hat\theta_t^{ls}.$

%\newpage

%%%%%%%%%%%%%%%%%%%%%%%%%%%%%%%%%%%%%%%%%%%%%%%%%%%%%%%%%%%%%%%%%%%%%%%%
%  7. Concluding remarks
%%%%%%%%%%%%%%%%%%%%%%%%%%%%%%%%%%%%%%%%%%%%%%%%%%%%%%%%%%%%%%%
\section {Concluding remarks}

This is a final part of a series of three papers (see Sharia (2006a) and Sharia (2006b)).
 We have  introduced
 estimation procedures  \eqref{rec2} which are recursive in the sense that each successive
 estimator is obtained from the previous one by a simple adjustment.
%An estimator of this class has a simple recursive
%structure so that there is no need of storing all the data
% and  repeating  calculations  each time a new observation
% is acquired.
To guarantee the convergence one has to impose global restrictions
on the functions in \eqref{rec2} (w.r.t. the parameter $\theta$)
such as a monotonicity type assumption and a restriction on the
growth at infinity (see Sharia (2006a)). This is the price one has to pay for the nice
recursive structure. Once the convergence is ensured,  the rate of convergence
(see Sharia (2006b)) and asymptotic
linearity can be deduced from local (in $\theta$) conditions.
Also, results presented
give an explicit way of constructing  a normalising sequence to ensure
local asymptotic linearity. The
rest relies on the ergodicity  of the model. Asymptotic properties
such as asymptotic distribution and efficiency of recursive (as
well as non-recursive) estimators depend on limit theorems
possessed by the model. For example,  in the i.i.d. case (see
Corollary 4.1), the central limit theorem and the law of large
numbers  imply
that the corresponding recursive procedures are asymptotically
normal and, in addition, the likelihood procedure is
asymptotically efficient. In general, one can obtain asymptotic
distribution and efficiency from asymptotic linearity (Lemma
3.1) and an appropriate central limit theorem.

The model considered in the paper is very general as we do not
impose any  preliminary  restrictions on probabilistic nature of
the observation process and cover a wide class of nonlinear
recursive procedures for estimation of a  multidimensional parameter.
The results are new even for the case of a scalar parameter and provide
a new insight even for the case of i.i.d. observations.

 While the advantage of this approach is its
universality, verification of the conditions may be a nontrivial
matter in some models. Examples considered give a flavour of what
is usually involved in this process and show where our
restrictions come from.  It is worth mentioning,
that even in the cases where one has difficulties with verifying
our conditions, the results of the paper can be used to determine
the form of  a recursive procedure (in fact, an algorithm, see
Remark 3.2 (iv)--(vi)), which is expected to have the same asymptotic
properties as the corresponding non-recursive one defined as a
solution of the equation \eqref{esteqg}.

\bigskip

\bigskip
\newpage

\begin{center}
{\large \bf APPENDIX A}
\end{center}

%\appendix
%\section{Appendix}

\noindent
{\bf Proof of Proposition 3.1}
To simplify notation we drop the fixed  argument or  the index
 $\theta$ in some of the expressions below.
\medskip

\noindent
To prove {\bf (a)},
 denote
$$
\chi_s=A_s[\tr\Gm_s(\theta)\Dl_{s-1}+R_s(\theta,\Dl_{s-1})]
$$
and
$$
{\cal G}_t=A_t^{-1}\sum_{s=1}^t\;\;
[\tr\Gm_s(\theta)\Dl_{s-1}+R_s(\theta,\Dl_{s-1})]=
A_t^{-1}\sum_{s=1}^t\;\;A_s^{-1}\chi_s.
$$
Applying the formula (summation by parts)
%                    PartSum
$$
\sum_{s=1}^t D_s\Dl C_s=D_tC_t - \sum_{s=1}^t \Dl D_s C_{s-1}, \;\;
C_0=0=D_0,
%\leqno{(A3)}
$$
 with
$C_s=\sum_{m=1}^s \chi_m$ and
$D_s=A_s^{-1} $ we obtain
$$
{\cal G}_t=
A_t^{-2}\sum_{s=1}^t
\chi_s
-A_t^{-1}\sum_{s=1}^t
\tr A_s^{-1}\sum_{m=1}^{s-1}
\chi_m.
$$
Then,
$
\tr A_s^{-1}=A_s^{-1}-A_{s-1}^{-1}=-A_s^{-1}(A_s-A_{s-1})A_{s-1}^{-1}
=-\tr A_s A_s^{-1}A_{s-1}^{-1},
$
where  the last equality follows since $A_s$ is diagonal.
Therefore,
$$
{\cal G}_t=
A_t^{-2}\sum_{s=1}^t
\chi_s
+A_t^{-1}\sum_{s=1}^t
\tr A_s \left\{ A_s^{-1}
A_{s-1}^{-1}\sum_{m=1}^{s-1}
\chi_m\right\}.
$$
Finally, since  $A_t$'s are diagonal with non-decreasing elements,
applying  the Toeplits Lemma to the components of the right hand side of latter formula we obtain that
${\cal G}_t \to 0.$
%%%%%%%%%%%%%

\medskip
\noindent
To prove {\bf (b)} and {\bf (c)}
denote ~ $M_t:=\sum_{s=1}^t {\cal E}_s.$
 Since $\bf \psi \in \Psi^M,$ it follows from that $M_t$ is a martingale.
Denote by $M^{(j)}_t$ the $j$-th component of  $M_t.$
Then  the square characteristic  ${\langle M^{(j)}\rangle}_t$ of
the martingale $M^{(j)}_t$ is
$$
{\langle M^{(j)}\rangle}_t=\sum_{s=1}^t
 E_\theta \left\{\left(
{\cal E}_s^{(j)}\right)^2 \mid{\cf}_{s-1}\right\}
$$
and, by  (LL2), ~
$
 \sum_{s=1}^\infty  {\tr {\langle M^{(j)}\rangle}_s} /{(A_s^{(jj)})^2} < \infty.
$
 It therefore follows that
$M^{(j)}_t/{A_t}^{(jj)}\to 0$  $P^\theta$ -a.s. (see e.g., Shiryayev (1984), Ch.VII, \S 5, Theorem 4).
This proves (c).  Now, use of the
Lenglart-Rebolledo inequality (see, e.g., Liptser and Shiryayev
(1989), Ch.1, $\S$9) yields
$$
P^\theta\left\{ (M^{(j)}_t)^2 \ge K^2 \left(
{A_t}^{(jj)}\right)^{2} \right\}\le \frac{\ve}{K}+
P^\theta\left\{ {\langle M^{(j)}\rangle}_t  \;\; \ge \; \ve \left(
{A_t}^{(jj)}\right)^{2} \right\}
$$
for each  $K>0$ and $\ve>0.$
Then, by (L2), ~
 $
 {\langle M^{(j)}\rangle_t}/( {A_t}^{(jj)})^2 \;
\to  \; 0 $ in probability $P^\theta$.
This implies that
$M^{(j)}_t/{A_t}^{(jj)}\to 0$ in probability $P^\theta  $ and so,
since $A_t$ is diagonal, (2) follows.
$\diamondsuit$
%%%%%%%%%%%%%%%%%%%%%%%%%%%%%%%%%%%%%%%%%%%%%%%%%%%%%%%%%%%%%%%%%%%%%%%%

\bigskip

%%%%%%%%%%%%%%%%%%%%%%%%%%%%%%%%%%%%%%%%%%%%%%%%%5555555555555555
\noindent {\bf Proof of Corollary 4.1} Using Corollary 4.1 in Sharia (2006a)
it follows that  (I) and
 (II) imply   $(\hat  \theta_t-\theta) \to
 0$.
We have $ \Gm_t(\theta)= t\gm(\theta) $ and $b(\theta, u)=\int
\p(\theta+u,z)f(\theta,z)\mu(\,dz). $ It is easy to see that
 (II) implies  (B2) from Corollary 4.1 in Sharia (2006b), and
 (V) implies that
(B1) of the same Corollary  holds with $C_\theta={{\bf 1}}.$ So, for any $0<\dl<1/2,$
%
%%%%%%%%%%%%%%%%%%%%%%%%%%%%%%%%%    (A1)
$$
t^\dl(\hat  \theta_t-\theta) \to 0
\leqno{(A1)}
$$
Let us check that conditions of Lemma 3.1 are also satisfied
with $ A_t={\sqrt t} {{\bf 1}}. $ Condition (EE) trivially holds.
 According to Proposition 3.1, condition (1) follows from
 (L1). To check (L1), it
is sufficient to show  that
%                             (A2)
$$
\frac 1t\sum_{s=1}^t[\gm(\theta)\Dl_{s-1}+R(\theta,\Dl_{s-1})]\sqrt
s \to 0,
\leqno{(A2)}
$$
where
$$
R(\theta, u)= R_t(\theta, u)=\gm(\theta)\gm^{-1}(\theta+u) \int
\p(\theta+u,z)f(\theta,z)\mu(\,dz).
$$
 By (V), ~
 $R(\theta, u)=-\gm(\theta)u+\gm(\theta)\gm^{-1}(\theta+u)\alpha^\theta(u)$ ~ and
$$
[\gm(\theta)\Dl_{s-1}+R(\theta,\Dl_{s-1})]\sqrt s
 =\sqrt s
\gm(\theta)\gm^{-1}(\theta+\Dl_{s-1})\alpha^\theta(\Dl_{s-1})
=\sqrt s\|\Dl_{s-1}\|^{1+\ve}\dl_s,
$$
where, by (III) and (V),
$\dl_s=\gm(\theta)\gm^{-1}(\theta+\Dl_{s-1})
\alpha^\theta(\Dl_{s-1})/\|\Dl_{s-1}\|^{1+\ve}\to 0$. Then,
$$
\sqrt s\|\Dl_{s-1}\|^{1+\ve}\dl_s=\sqrt{\frac {s}{s-1}}
\left((s-1)^{\frac1{2(1+\ve)}}\|\Dl_{s-1}\|\right)^{1+\ve}\dl_s
$$
which, by (A1) (since $1/({2(1+\ve)}) < 1/2$) converges to zero. Therefore,
(A2) is now a consequence of the Toeplits Lemma.

For the process ${\cal E}_s(\theta)$ from (L2) (since
$\|u-v\|^2\le 2\|u\|^2+2\|v\|^2$), we have
$$
\|{\cal E}_s(\theta)\|^2=\|\gm(\theta) \gm^{-1}(\theta+\Dl_{s-1})
\left(\p(\theta+\Dl_{s-1},X_s)-b(\theta,\Dl_{s-1})\right)
-\p(\theta,X_s)\|^2
$$
$$
\le 2\|\gm(\theta) \gm^{-1}(\theta+\Dl_{s-1})
\p(\theta+\Dl_{s-1},X_s) -\p(\theta,X_s)\|^2 +2\|\gm(\theta)
\gm^{-1}(\theta+\Dl_{s-1}) b(\theta,\Dl_{s-1})\|^2.
$$
From (III) and (V) we obtain that $\left(\gm(\theta)
\gm^{-1}(\theta+\Dl_{s-1})-{{\bf 1}}\right)\to 0$ ~ and
$b(\theta,\Dl_{s-1})\to 0$ as $s\to \infty$. So, using (IV), it is easy to see that$
 E_\theta \left\{ \left( {\cal E}^{(j)} _s(\theta)\right)^2
\mid{\cf}_{s-1}\right\} \to 0. $ Since
$
{(A_t^{(jj)}(\theta))^2}
 = t,$
(L2) follows from  the Toeplitz lemma.

Therefore, the conditions of Lemma 3.1 hold for
$A_t(\theta)=\sqrt{t}.$ This implies that  $\sqrt{t}(\hat \theta_t-\theta_t^*)\to 0$
in probability $P^\theta,$ where
$$
\theta_t^*=\frac1{t\gamma(\theta)}\sum_{s=1}^t {\psi}_s(\theta, X_s).
$$
  The asymptotic normality now obviously follows
from the central limit theorem for i.i.d.  random variables.
$\diamondsuit$

%%%%%%%%%%%%%%%%%%%%%%%%%%%%%%%%%%%%%%%%%%%%%%%%%%5555555555555555

%%%%%%%%%%%%%%%%%%%%%%%%%%%%%%%%%%%%%%%%%%%%%%%%%%%%%%%%%%%%%%%%%%%%%%%%%%%
\bigskip

\noindent
{\bf Corollary A1}
{\it
Suppose that  $H_t\to \infty$ and $\hat\theta_t$ is derived by
\eqref{RecAdd}. Denote  $\Dl_t=\hat\theta_t-\theta$,
$l_t(\theta)=m(X_t,X_{t-1})-\dot\gamma(\theta)h(X_{t-1})$, and
suppose also that
\begin{description}
\item[(I)]
$$
H_t^{-1/2}\sum_{s=1}^t {\cal E}_s(\theta) \to  0,
$$
where
$$
{\cal E}_s(\theta) =\frac{\ddot\gamma(\theta+\Dl_{s-1})-\ddot\gamma(\theta)}{\ddot\gamma(\theta+\Dl_{s-1})}
l_s(\theta);
$$
\item[(II)] one of the following two conditions are satisfied;
$$
H_t^{-1/2}\sum_{s=1}^t \tr H_s{\cal C}_s(\theta) \to  0,
$$
OR
$$
H_t^{-1}\sum_{s=1}^t \tr H_s ~ H_s^{1/2}{\cal C}_s(\theta) \to  0,
$$
where
$$
{\cal C}_s(\theta)=\frac{\ddot\gamma(\theta+\Dl_{s-1})-\ddot\gamma(\theta+\tilde \Dl_{s-1})}
{\ddot\gamma(\theta+\Dl_{s-1})} \Dl_{s-1}
$$
and $\tilde\Dl_t$ is a predictable process with  $|\tilde\Dl_t|\le |\Dl_t|.$
\end{description}
Then \eqref{aslinexp} holds, i.e.,  the estimator $\hat\theta_t$ is asymptotically linear.
}
\medskip
%%%%%%%%%%%%%%%%%%%%%%%%%%%%%%%%%%%%%%%%%%%%%%%%%%%%%%%%%%%%%%%%%%%%%%%%%%

\noindent
 {\bf Proof.}
 Let us check the conditions  of Lemma 3.1 for  $\psi_t(\theta)=l_t(\theta),$
%%%%%%%%%%%%%%%%%%%%%%%%%%%%%%%%%%%%%%%%%%%%%%%%%%%%%%%%%%%%%%%%%%%%%%%%%%%
%                               A2
%%%%%%%%%%%%%%%%%%%%%%%%%%%%%%%%%%%%%%%%%%%%%%%%%%%%%%%%%%%%%%%%%%%%%%%%%
$$
\Gamma_t(\theta)=I_t(\theta)=\ddot\gamma(\theta)H_t
\leqno{(A3)}
$$
and  $A_t(\theta)=H_t^{1/2}.$
 Since $l_t(\theta)$
is a martingale-difference, we have
$ E_\theta \left\{m(X_t,X_{t-1})\mid{\cf}_{t-1}\right\}
=\dot\gamma(\theta)h(X_{t-1})$ and so
%%%%%%%%%%%%%%%%%%%%%%%%%%%%%%%%%%%%%%%%%%%%%%%%%%%%%%%%%%%%%%%%%%%%%%%%%%%
%                               A2
%%%%%%%%%%%%%%%%%%%%%%%%%%%%%%%%%%%%%%%%%%%%%%%%%%%%%%%%%%%%%%%%%%%%%%%%%
$$
b_t(\theta,u)=E_\theta\left\{l_t(\theta+u)\mid
{\cal{F}}_{t-1}\right\}=
h(X_{t-1})\left(\dot\gamma(\theta)-\dot\gamma(\theta+u)\right)
\leqno{(A4)}
$$
and
$$
R_t(\theta,u)=
\frac{\ddot\gamma(\theta)}{\ddot\gamma(\theta+u)}h(X_{t-1})({\dot\gamma(\theta)-\dot\gamma(\theta+u)})
=-\frac{\ddot\gamma(\theta)}{\ddot\gamma(\theta+u)}h(X_{t-1})\ddot\gamma(\theta+\tilde u)u
$$
where $|\tilde u|\le |u|.$
Then, since $\tr\Gamma_t(\theta)=\tr I_t(\theta)=h(X_{t-1})\ddot\gamma(\theta)$ we have
$$
\tr\Gamma_t(\theta)u+R_t(\theta,u)
=h(X_{t-1})\ddot\gamma(\theta)
\frac{\ddot\gamma(\theta+u)-\ddot\gamma(\theta+\tilde u)}{\ddot\gamma(\theta+u)}u.
%\le\tr\Gamma_t(\theta)\frac{K_\theta |u|^\ve}{\ddot\gamma(\theta+u)}.
$$
Now, since $\tr H_t=h(X_{t-1}),$
it is easy to see that the first condition in   (II) implies  (1) in Lemma 3.1 and
the second condition in (II) implies  (L1) in Proposition 3.1. Therefore, (1) holds.

To verify (2), consider  the process ${\cal E}_s(\theta)$ defined in
(2).  Using (A3) and (A4), it is easy to see that
$$
{\cal E}_s(\theta)
=\left(1-\frac{\ddot\gamma(\theta)}{\ddot\gamma(\theta+\Dl_{s-1})}\right)
(m(X_s,X_{s-1})-\dot\gamma(\theta)h(X_{s-1}))
$$
$$
=\frac{\ddot\gamma(\theta+\Dl_{s-1})-\ddot\gamma(\theta)}{\ddot\gamma(\theta+\Dl_{s-1})}
l_s(\theta). \leqno{(A5)}
$$
This shows that (I) implies (2).
$\diamondsuit$

\bigskip

%\medskip
%%%%%%%%%%%%%%%%%%%%%%%%%%%%%%%%%%%%%%%%%%%%%%%%%%%%%%%%%%%%%%%%%%%%%%%%%%%%%%%%
%
%%%%%%%%%%%%%%%%%%%%%%%%%%%%%%%%%%%%%%%%%%%%%%%%%%%%%%%%%%%%%%%%%%%%%%%%%%
\noindent
 {\bf Proof of Proposition 4.3}
Since, by (iii),  $\ddot\gamma(\cdot)$ is obviously a continuous   function,
 condition (M2) of Proposition 4.1 in Sharia (2006b) holds. Also, (M1) in the same proposition
 obviously follows from (i).   So, it follows that all the conditions
 of Proposition 4.1 and Corollary 4.2 in Sharia (2006b) are satisfied implying that
 $H_t^\delta(\hat\theta_t - \theta )\to 0\;\; (P^\theta$-a.s.). Also, by (i),
 $\tr H_t/H_{t-1}={h(X_{t-1})}/{H_{t-1}}\to 0$  implying that
$ H_t/H_{t-1}=1+\tr H_t/H_{t-1} \to 1.$
So,
$$
H_t^\delta \Dl_{t-1}=H_t^\delta(\hat \theta_{t-1} - \theta)\to 0.
\leqno{(A6)}
$$
To establish asymptotic linearity, let us verify the conditions of Corollary A1 is satisfied.
Since $\Dl_{s-1}=\hat \theta_{s-1} -\theta \to 0\;\; (P^\theta$-a.s.) and
$|\tilde \Dl_{s-1}|\le |\Dl_{s-1}|$,  by (iii) we obtain that
$|\ddot\gamma(\theta+\Dl_{s-1})-\ddot\gamma(\theta+\tilde \Dl_{s-1})|\le 2K_\theta |\Dl_{s-1}|^{\ve_\theta}$
eventually.
So,
$$
|H_s^{\frac12}{\cal{C}}_s(\theta)|=H_s^{\frac12}
\frac{|\ddot\gamma(\theta+\Dl_{s-1})-\ddot\gamma(\theta+\tilde \Dl_{s-1})|
 |\Dl_{s-1}|}
{\ddot\gamma(\theta+\Dl_{s-1})}
\le \frac{2K_\theta H_s^{\frac12}|\Dl_{s-1}|^{1+{\ve_\theta}}}{\ddot\gamma(\theta+\Dl_{s-1})}
$$
eventually.
Now,
$$
H_s^{\frac12}|\Dl_{s-1}|^{1+{\ve_\theta}}=
|H_s^{\frac1{2(1+{\ve_\theta})}}(\hat\theta_{s-1}-\theta)|^{1+{\ve_\theta}}\to 0,
$$
by (A6) since  $\frac1{2(1+{\ve_\theta})}< \frac12.$
 So, since $\ddot \gamma(\cdot)$ we obtain that
 $|H_s^{\frac12}{\cal{C}}_s(\theta)|\to 0.$ Therefore,  by the Toeplits Lemma, the second condition of (II) holds.

Now, since ${\cal E}_s(\theta)$ is a martingale-difference,
to verify (I), it is sufficient to show that (see e.g., Shiryayev (1984), Ch.VII, \S 5, Theorem 4)
$$
\sum_{s=1}^\infty\frac{E_\theta \left\{ {\cal E}_s^2(\theta)
\mid{\cf}_{s-1}\right\}}{H_s} < \infty.
$$
Since $
E_{\theta}\{l^{2}_s(\theta)\mid
{\cf}_{s-1}\}=\ddot\gamma(\theta)h(X_{s-1})=\ddot\gamma(\theta)\tr H_s,
$
the above series  can be rewritten as
$$
\sum_{s=1}^\infty \frac{\tr H_s}{H_s} \ddot\gamma(\theta)
\left(\frac{\ddot\gamma(\theta+\Dl_{s-1})-\ddot\gamma(\theta)}{\ddot\gamma(\theta+\Dl_{s-1})}\right)^2
=\ddot\gamma(\theta)
\sum_{s=1}^\infty \frac{\tr H_s}{H_s^{1+{\ve_\theta}/2}}r_s
$$
where, by (iii),
$$
r_s=
\frac{\left(\ddot\gamma(\theta+\Dl_{s-1})-\ddot\gamma(\theta)\right)^2H_s^{{\ve_\theta}/2}}
{\ddot\gamma^2(\theta+\Dl_{s-1})}
\le K_\theta^2
\frac{|\Dl_{s-1}|^{2{\ve_\theta}}H_s^{{\ve_\theta}/2}}
{\ddot\gamma^2(\theta+\Dl_{s-1})}
=K_\theta^2
\frac{(|\Dl_{s-1}|H_s^{1/4})^{2{\ve_\theta}}}
{\ddot\gamma^2(\theta+\Dl_{s-1})}.
$$
Now, using   (A6) and continuity of $\ddot\gamma(\cdot)$ we deduce that $r_s \to 0.$
Also,
$$
\sum_{s=1}^\infty \frac{\tr H_s}{H_s^{1+{\ve_\theta}/2}} < \infty
$$
(see Sharia (2006b),
Appendix A, Proposition A2),
implying   that the above series converge which completes the proof.  $\diamondsuit $

\bigskip
%%%%%%%%%%%%%%%%%%%%%%%%%%%%%%%%%%%%%%%%%%%%%%%%%%%%%%%%%%%%%%%%%%%
%\section                  {R E F E R E N C E S }
%%%%%%%%%%%%%%%%%%%%%%%%%%%%%%%%%%%%%%%%%%%%%%%%%%%%%%%%%%%%%%%%%%%
\begin{center}
       {\Large{\bf  REFERENCES}}
\end{center}
%\vskip+0.5cm
%\begin{enumerate}
%\begin{itemize]
\begin{description}
\item
 {Barndorff-Nielsen,} O.E. and {Sorensen,} M. (1994).
 A review of some aspects of asymptotic likelihood theory for
 stochastic processes. {\it International Statistical Review.}
 {\bf 62,} 1, 133-165.
\item
     {Basawa,} I.V. and  {Scott,} D.J. (1983).
     {\it Asymptotic Optimal Inference for  Non-ergodic Models.}
      Springer-Verlag, New York.

\item
   {Campbell,} K.(1982). Recursive computation of
     M-estimates for the parameters of a finite autoregressive
     process, {\it Ann. Statist.}, {\bf 10}, 442-453.
 \item
 { Englund,} J.-E.,   {Holst,} U., and {Ruppert,} D.(1989).
  Recursive estimators for
         stationary, strong mixing processes -- a representation
         theorem and asymptotic distributions,
         {\it Stochastic Processes Appl.},  {\bf 31},  203--222.

\item
{Fabian,} V. (1978).  On asymptotically  efficient
         re\-cur\-sive es\-ti\-ma\-tion.
         {\it Ann. Statist.} {\bf 6}, 854-867.
\item
Feigin, P.D. (1985). Stable convergence for semimartingales. {\it Stoch.  Proc. Appl.}
{\bf 19}, 125--134.
\item
 {Hall,} P. and  {Heyde,} C.C. (1980). {\it Martingale Limit Theory and Its
          Application.} Academic Press, New York.
 \item
  {Hampel,} F.R.,   {Ronchetti,} E.M.,   {Rousseeuw,}
  P.J., and {\sc Stahel,} W. (1986).
         {\it Robust Statistics - The Approach Based on Influence
         Functions}, Wiley, New York.
\item
{Heyde,} C.C. (1997). {\it
Quasi-Likelihood and Its Application: A General Approach to
Optimal Parameter estimation.} Springer-Verlag, New York.
 \item
Hutton, J.E. and Nelson, P.I. (1986). Quasi-likelihood estimation
for semimartingales. {\it Stoch.  Proc. Appl.} {\bf 22}, 245--257.
\item
Jacod, J. and {Shiryayev,} A.N. (1987). {\it Limit Theorems
for Stochastic Processes.} Heidelberg, Springer.
\item
{Jure$\check{\mbox{{c}}}$kov$\acute{\mbox{{a}}}$,} J.  and
    {Sen,} P.K.   (1996). {\it Robust Statistical Procedures -
    Asymptotics and Interrelations,} Wiley, New York.

\item
  {Khas'minskii,} R.Z. and {Nevelson,} M.B.     (1972).
  {\it Stochastic Approximation and
          Recursive Estimation.} Nauka, Moscow.
 \item
    {Lazrieva,}  N., {Sharia,} T. and  {Toronjadze,} T.(1997).
     The Robbins-Monro type
     stochastic differential equations. I. Convergence of solutions,
      {\it Stochastics and Stochastic Reports,} {\bf 61}, 67--87.
\item
    {Lazrieva,}  N., {Sharia,} T. and  {Toronjadze,} T.(2003).
     The Robbins-Monro type
     stochastic differential equations. II.  Asymptotic behaviour of solutions,
      {\it Stochastics and Stochastic Reports,} {\bf 75}, 153--180.
\item
    {Lazrieva,} N. and  {Toronjadze,} T.  (1987). Ito-Ventzel's formula for
semimartingales, asymptotic properties of MLE and recursive
    estimation,  {\it Lect. Notes in Control and Inform. Sciences, 96, Stochast.
    diff.  systems, H.J, Engelbert, W. Schmidt (Eds.),} (pp. 346--355). Springer.

  \item
 {Ljung,} L.  {Pflug,} G. and {Walk,} H. (1992).
 {\it Stochastic Approximation and
         Optimization of Random Systems,}  Birkh\"auser, Basel.
\item
  {Ljung,} L. and {Soderstrom,} T. (1987). {\it Theory and
Practice of Recursive Identification,} MIT Press.

\item
        {Prakasa Rao,} B.L.S. (1999). {\it
        Semimartingales and their Statistical Inference.}
     Chapman $\&$ Hall,   New York.
     \item
     {Rieder,} H. (1994). {\it Robust Asymptotic Statistics,}
         Springer--Verlag,  New York.

\item
    {Sharia,} T. (1998). On the recursive parameter estimation for the
          general discrete time statistical model. {\it
          Stochastic Processes Appl.}
          {\bf 73}, {\bf 2}, 151--172.
 \item
{Sharia,} T.   (2006a). Recursive  parameter estimation: Convergence.
{\it Statistical Inference for Stochastic Processes} (in press).\\
(see also {\em http://personal.rhul.ac.uk/UkAH/113/ConvA.pdf}).

\item
 {Sharia,} T.   (2006b).
  Rate of  convergence in recursive parameter estimation procedures (submitted).
  ({\em http://personal.rhul.ac.uk/UkAH/113/GmjA.pdf}).
  \item
  {Shiryayev,} A.N. (1984). {\it Probability.}
           Springer-Verlag, New York.
\end{description}
%\end{enumerate}
%\end{itemize}

%
\end{document}